\documentclass[12pt]{amsart}            
\usepackage{amssymb}
\usepackage{amscd}
\usepackage{pstricks}
\usepackage{pst-node}
\newtheorem{thm}{Theorem}[section]
\newtheorem{cor}[thm]{Corollary}
\newtheorem{lem}[thm]{Lemma}

\newtheorem{prop}[thm]{Proposition}
\theoremstyle{definition}

\newtheorem{example}[thm]{Example}

\newtheorem{rem}[thm]{Remark}

\theoremstyle{remark}

\numberwithin{equation}{section}

\newcommand{\trans}[1]{{}^t\kern-.2em{#1}}
\newcommand{\ytrans}[1]{{}^t\kern-.11em{#1}}
\newcommand{\Trans}[1]{{}^T\kern-.2em{#1}}
\newcommand{\lsup}[2]{{}^{#1}\kern-.1em{#2}}

\newcommand{\M}{\mathbf{\M}}

\renewcommand{\tilde}[1]{\widetilde{#1}}

\DeclareFixedFont{\bgn}{OT1}{cmr}{m}{n}{20.74}
\DeclareFixedFont{\bgi}{OT1}{cmr}{m}{it}{20.74}

\newcommand{\bigzerou}{\smash{\lower1.7ex\hbox{\bgi O}}}

\def\eqnarray{%
   \stepcounter{equation}%
   \def\@currentlabel{\p@equation\theequation}%
   \global\@eqnswtrue
   \m@th
   \global\@eqcnt\z@
   \tabskip\@centering
   \let\\\@eqncr
   $$\everycr{}\halign to\displaywidth\bgroup
       \hskip\@centering$\displaystyle\tabskip\z@skip{##}$\@eqnsel
      &\global\@eqcnt\@ne \hfil$\displaystyle{{}##{}}$\hfil
      &\global\@eqcnt\tw@ $\displaystyle{##}$\hfil\tabskip\@centering
      &\global\@eqcnt\thr@@ \hb@xt@\z@\bgroup\hss##\egroup
        \tabskip\z@skip
      \cr}
\makeatother
\makeatletter
\def\varin{\mathrel{\mathpalette\@varin\relax}}
\def\@varin#1{%
   \hbox{\setbox\z@\hbox{\m@th$#1\cup$}%
       \def\reserved@a{bold}%
       \dimen@\ifx\reserved@a\math@version .3\else .2\fi\p@
       \kern.5\wd\z@\kern-\dimen@
       \vrule\@width2\dimen@\@height1.08\ht\z@\@depth\z@
       \kern-\dimen@\kern-.5\wd\z@
       \box\z@}}
\makeatother

\begingroup
    \divide\count0 by 60
    \multiply\count0 by 60
    \advance\count1 by -\count0
    \ifnum\count2>11
         \ifnum\count2>12 \advance\count2 by -12\fi
         \def\ampm{\,pm}%
    \else
         \ifnum\count2=0 \advance\count2 by 12\fi
         \def\ampm{\,am}%
    \fi
    \xdef\daytime{%
         \ifnum\count2<10 0\fi \the\count2:%
         \ifnum\count1<10 0\fi \the\count1
         \ampm
    }%
    \xdef\Daytime{
         \the\count2:
         \ifnum\count1<10 0\fi \the\count1   
         \ampm                               
    }
\endgroup

\begin{document}
 
\keywords{Heisenberg group, zeta-regularized determinant,
Laplacian, heat kernel, Kroneker's second limit formula, modified
Bessel function, Poisson's summation formula}

\title[Determinant of Laplacians] 
{Determinant of Laplacians on Heisenberg manifolds}

\author{Kenro Furutani and Serge de Gosson}
\address{Kenro Furutani \endgraf 
Department of Mathematics \endgraf 
Faculty of Science and Technology \endgraf 
Science University of Tokyo \endgraf 
2641 Noda, Chiba (278-8510)\endgraf 
Japan \endgraf}
\email{furutani{\_}kenro@ma.noda.tus.ac.jp}
\address{Serge de Gosson \endgraf      
Department of Mathematics \endgraf 
V\"axj\"o University \endgraf
SE-351 95 V\"axj\"o \endgraf 
Sweden\endgraf}
\email{sdg@vax.se}



\bigskip
\bigskip
\maketitle
\begin{abstract}
We give an integral representation of the
{\it zeta-regularized determinant} of Laplacians on
three dimensional Heisenberg manifolds, and
study a behavior of the values when we deform
the uniform discrete subgroups. Heisenberg manifolds are
the total space of a fiber bundle with a torus as the base space
and a circle as a typical 
fiber, then the deformation of the uniform discrete
subgroups means that the ``radius'' of the fiber goes to zero.
We explain the lines of the calculations 
precisely for three dimensional
cases and state the corresponding results 
for five dimensional Heisenberg manifolds. 
We see that the values themselves are of the
product form with a factor which is that of the flat torus. 
So in the last half of this paper we derive general 
formulas of the zeta-regularized determinant
for product type manifolds of two Riemannian manifolds, 
discuss the formulas for flat tori and explain
a relation of the formula for the two dimensional 
flat torus and the {\it Kroneker's second limit formula}.
\end{abstract}

\maketitle
\tableofcontents
\thispagestyle{empty}

\section*{Introduction}
In the conformal field theory, especially in the 
string theory an infinite dimensional analogue of the
determinant for elliptic operators appears and it is 
considered in the framework of the analytic continuation 
method of the spectral zeta-function of elliptic operators. 
In the string theory it seems to be most interesting to
calculate the explicit values for compact Riemannian 
surfaces and also from the mathematical point of view 
the quantity relates with various aspects
of special functions already in such two dimensional cases.
As for the values for the spheres there are relations 
with the values of the Riemann $\zeta$-functions at the negative 
integers, for the two dimensional flat torus it is expressed 
by using the famous formula, 
so called ``{\it Kroneker's  second limit formula}'', 
and for the cases of compact surfaces with constant 
negative curvature it was calculated by using Selberg's 
trace formula and deep properties of modified Bessel functions
(\cite{Mc}, \cite{Va}, \cite{DP}, \cite{BS}, \cite{Fo}). 
Thus for all two dimensional cases 
we already know the values of the zeta-regularized determinant
of the Laplacians or relations with other quantities. 
A purpose of this note is to give an integral representation
of the zeta-regularized determinant 
for three dimensional Heisenberg manifolds 
and state the corresponding results for five dimensional
Heisenberg manifolds.  
It will be possible to give similar expressions
for higher dimensional cases, 
but we restrict ourselves to these two cases
and also we restrict ourselves 
to deal with a certain kind of uniform discrete
subgroups of the Heisenberg group, since even in these cases they
contain
all necessary features for determining the values for any cases 
and make us
the calculations to be simple. Then in the last half of this paper
we derive a general formula of the zeta-regularized determinant
for manifolds of the product form of two Riemannian manifolds
and the formulas for flat tori of two, three and four dimensions.

In $\S 1$ we gather up the basic data of the spectrum of the three
dimensional Heisenberg group and Heisenberg manifolds. In $\S 2$
first, we explain the zeta-reguralized determinant of the 
Laplacian and give a calculation for the three 
dimensional Heisenberg manifolds based on an integral
representation of the spectral zeta-function 
(= the Mellin transformation of the trace of heat kernel
divided by a Gamma function). 
In $\S 3$ as an application
of an integral representation of the spectral zeta-function given in $\S 2$
we give expressions of the all coefficients of the asymptotic
expansion of the heat kernel for the three dimensional Heisenberg manifolds.
In $\S 4$ we state the corresponding results in $\S 2$ and $\S 3$ 
for five dimensional Heisenberg manifolds. 
In $\S 5$ we give a general expression
of the zeta-regularized determinant of the Laplacian
on the product of two Riemannian manifolds.
In $\S 6$ we give a precise
form of the formula derived in $\S 5$ for a product type manifold with $S^1$.
Finally
in $\S 7$ we give such formulas for two, three and four dimensional
flat tori and explain a relation of the formula for
two dimensional flat torus and 
the {\it Kroneker's second limit formula}.
\section{Spectrum of three dimensional Heisenberg manifolds}

Let $H_{3}$ 
be the three dimensional Heisenberg group:
\begin{equation}
H_{3}=\left\{ g=g\left(x,y,z\right)=
\begin{pmatrix}
1 & x & z\\
0 & 1 & y\\
0 & 0 & 1
\end{pmatrix}
\,|\,x,y,z\in\mathbb{R}\right\}.
\end{equation}
The Lie algebra is  
\begin{equation}
\mathfrak{h}_3= \left\{ X=X\left(x,y,z\right)=
\begin{pmatrix}
0 & x & z\\
0 & 0 & y\\
0 & 0 & 0
\end{pmatrix}
\,|\,x,y,z\in\mathbb{R}\right\}.
\end{equation}
It is decomposed into a direct sum in the form of
\begin{equation}
\mathfrak{h}_3=\mathfrak{g}_+ \oplus\mathfrak{g}_- \oplus\mathfrak{z}, 
\end{equation}
where 
$$\mathfrak{g}_+=\left\{\begin{pmatrix}
0 & x & 0\\
0 & 0 & 0\\
0 & 0 & 0
\end{pmatrix}
\,|\,x\in\mathbb{R}\right\},
$$  
$$\mathfrak{g}_-=\left\{\begin{pmatrix}
0 & 0 & 0\\
0 & 0 & y\\
0 & 0 & 0
\end{pmatrix}
\,|\,y\in\mathbb{R}\right\}
$$ 
and  
$$\mathfrak{z}=\text{the center}=\left\{\begin{pmatrix}
0 & 0 & z\\
0 & 0 & 0\\
0 & 0 & 0
\end{pmatrix}
\,|\,z\in\mathbb{R}\right\}.
$$

We identify $\mathfrak{h}_3$ and $H_3$ through the exponential map.
Then the group multiplication of two elements 
$X =X\left(x,y,z\right)\in \mathfrak{h}_3$ and
$\tilde{X}
=\tilde{X}\left(\tilde{x},\tilde{y},\tilde{z}\right)\in\mathfrak{h}_3$
is given by
\begin{equation}
Y=X*\tilde{X}
= X+\tilde{X} +\frac{1}{2}[X,\tilde{X}],\,\text{that is},\,
\exp \left(Y\right)= \exp{X}\cdot \exp{\tilde X}.
\end{equation}

Left invariant Riemannian metrics are determined by its restriction
to the tangent space at the identity element ($\cong$ its Lie
algebra), and  among the left invariant Riemannian metrics 
we only consider such a metric that 
$X_1=\begin{pmatrix}
0 & 1 & 0\\
0 & 0 & 0\\
0 & 0 & 0
\end{pmatrix}$, 
$Y_1=\begin{pmatrix}
0 & 0 & 0\\
0 & 0 & 1\\
0 & 0 & 0
\end{pmatrix}$ 
and $Z_0 =z_0\cdot Z=z_0\cdot 
\begin{pmatrix}
0 & 0 & 1\\
0 & 0 & 0\\
0 & 0 & 0
\end{pmatrix}$ 
are being an orthonormal basis.  Also 
we only consider
uniform discrete subgroups 
$\Gamma_{\ell}$ of the following form
\begin{equation}
\Gamma_{\ell}=\left\{\begin{pmatrix}
1 & m & k/2{\ell}\\0 & 1 & n\\0 & 0 & 1\end{pmatrix}
\,|\,m,n,k\in\mathbb{Z}\right\},
\end{equation}
where $\ell$ is a positive integer.
These choices do not loose the essential features in treating with
the spectrum of the Laplacian on the quotient space,
$H_3/\Gamma_{\ell}$, so called, Heisenberg manifolds, 
since the spectrum is given more or less in a similar form
(\cite{GW}, we state them later).

Then the inverse image of $\Gamma_{\ell}$ by the exponential map is
\begin{equation}
\exp^{-1}\left(\Gamma_{\ell}\right)=\left\{
\begin{pmatrix}
0 & m & k/2{\ell}\\
0 & 0 & n\\
0 & 0 & 0 
\end{pmatrix}
\,|\,m,n,k\in\mathbb{Z}\right\},
\end{equation}
which is a direct sum of two uniform lattices
$\Gamma_B$ in $\mathfrak{g}_+\oplus\mathfrak{g}_-$
and $\Gamma_V\left({\ell}\right)$ in $\mathfrak{z}$ such that
\begin{equation}
\Gamma_B =\left\{\begin{pmatrix}
0 & m & 0\\
0 & 0 & n\\
0 & 0 & 0 
\end{pmatrix}
\,|\,m,n\in\mathbb{Z}\right\}
\end{equation}
and
\begin{equation}
\Gamma_V\left(\ell\right) =\left\{\begin{pmatrix}
0 & 0 & k/2{\ell}\\0 & 0 & 0\\0 & 0 & 0 \end{pmatrix}
\,|\,k\in\mathbb{Z}\right\}.
\end{equation}

The kernel $K_t\left(g,\tilde{g}\right)$ 
of $e^{-t\Delta}$, the heat kernel on the Heisenberg group 
$H_3$, is given by
\allowdisplaybreaks{\begin{align}
&K_t\left(g,\tilde{g}\right)
=K_t\left(x,y,z;\tilde{x},\tilde{y},\tilde{z}\right)\notag\\
&=\left(2\pi\right)^{-2}
\int_{-\infty}^{+\infty}e^{\sqrt{-1}\cdot\eta\cdot\left\{\tilde{z}-z+
\frac{1}{2}\left(\tilde{x}y-x\tilde{y}\right)\right\}}
\cdot e^{-t\left|\eta\right|^2}\times\\
&\times\frac{\left|\eta\right|}{2\sinh t\left|\eta\right|}
e^{
-\frac{\left|\eta\right|}{4}
\frac{\cosh t\left|\eta\right|}{\sinh t\left|\eta\right|}
\cdot\left\{\left(x-\tilde{x}\right)^2
+\left(y-\tilde{y}\right)^2\right\}}d\eta,\notag
\end{align}}
where we regard 
$\eta\in \mathfrak{z}^*=\mathbb{R}Z^*\cong\mathbb{R}$ and
$g=xX_1+yY_1+zZ$, $\tilde{g}=\tilde{x}X_1+\tilde{y}Y_1+\tilde{z}Z$
$\in$ $\mathfrak{h}_3 \cong H_3$
through the identification by the exponential map
$\exp:\mathfrak{h}_3 \to H_3$.

Then the heat kernel
$k_{H_3/\Gamma_{\ell}}\left(t;[g],[\tilde{g}]\right)$ 
on the Heisenberg manifold $H_3/\Gamma_{\ell}$
is expressed by making use of the heat kernel
$K_t\left(g,\tilde{g}\right)$
on the whole group, because of the invariance
$K_t\left(\gamma\cdot g,\gamma\cdot\tilde{g}\right)
=K_t\left(g,\tilde{g}\right), 
\,\gamma\in
H_3$:
\begin{equation}
k_{H_3/\Gamma_{\ell}}\left(t;[g],[\tilde{g}]\right)=
\sum\limits_{\gamma\in \Gamma_{\ell}} 
K_t\left(\gamma\cdot g,\tilde{g}\right).
\end{equation}

Its trace is calculated in the following form:
\begin{thm}{\bf (\cite{Fu})}
\allowdisplaybreaks{\begin{align*}
&\int_{H_3/\Gamma_{\ell}}
k_{H_3/\Gamma_{\ell}}\left(t;[g],[g]\right)dg
=\sum\limits_{\gamma\in\Gamma_{\ell}}\int_{F_{\ell}}
K_t\left(\gamma\cdot g, g\right)dg\\
&=\sum\limits_{\mu\in\Gamma_V\left(\ell\right)^{*},\mu\not=0}
\sum\limits_{m=0}^{\infty}
Vol\left(\mathfrak{g}_\oplus\mathfrak{g}_-/\Gamma_B\right)\Vert\mu\Vert
\cdot e^{-t\left\{4\pi^2{\Vert\mu\Vert}^2
+2\pi\left(2m+1\right)\Vert\mu\Vert\right\}}\\
&+\sum\limits_{\nu\in\Gamma_B^*}e^{-4\pi^{2}t{\Vert\nu\Vert}^{2}}.
\end{align*}}
\end{thm}
Here $F_{\ell}$ denotes a fundamental domain of the uniform
discrete subgroup $\Gamma_{\ell}$ and 
$\Gamma_B^*$ and $\Gamma_V\left(\ell\right)^*$ 
are dual lattices, i.e.,
\allowdisplaybreaks{\begin{align*}
\Gamma_B^*
&=\left\{ \nu\in\left(\mathfrak{g}_+\oplus\mathfrak{g}_-\right)^*\,|\,
\nu\left(\gamma\right)\in \mathbb{Z},\,
\text{for any}\,\gamma\in\Gamma_B\right\}\\
&\cong\{\nu=\nu_1 X_1^*+\nu_2 Y_1^*\,|\,\nu_i\mathbb{Z}\},
\,\, \,\left\{X_1^*,\,Y_1^*\right\}\,\,\text{are dual basis of}\,\,
\mathfrak{g}_+\oplus\mathfrak{g}_-,\\
\Gamma_V\left(\ell\right)^*
&=\left\{\mu\in\left(\mathfrak{z}\right)^*\,
\left|\right.\,\mu\left(\gamma\right)\in\mathbb{Z},\,\text{for any}\, 
\gamma\in\Gamma_V\left(\ell\right)\right\}\\
&\cong\left\{\mu=2\ell\cdot k Z_1^*\,|\,k\in\mathbb{Z}\right\},
\,\,\,\Vert\mu\Vert= 2\ell |k|. 
\end{align*}}

The spectrum of the Laplacian on $H_3/\Gamma_{\ell}$ is the union
of eigenvalues of two types ({\it a}) and ({\it b}):
\begin{prop}{\bf (\cite{GW})}
\allowdisplaybreaks{\begin{align*}
\text{(a)}:\quad & 4\pi^2{\Vert\mu\Vert}^2+2\pi\left(2m+1\right)
\Vert\mu\Vert,\,\,\mu=2\ell\cdot k Z_1^*\in \Gamma_V\left(\ell\right)^*\\
\qquad\qquad &\text{with the multiplicity equal to}\\
\qquad\qquad\qquad & 
Vol\left(\mathfrak{g}_+\oplus\mathfrak{g}_-/\Gamma_B\right)
\cdot \Vert\mu\Vert
= 1\cdot 2\ell\cdot |k| \,\left(k\in\mathbb{Z}\right)\\
\text{(b)}:\quad & 4\pi^2 \Vert \nu \Vert^2
=4\pi^2\left({\nu_1}^2+{\nu_2}^2\right),\,\, 
\nu=\nu_1 X_1^*+\nu_2 Y_1^*\in \Gamma_B^*.
\end{align*}}
\end{prop}

{}From the exact sequence
$$
0
\longrightarrow 
\left\{\begin{pmatrix}
1 & 0 & z\\
0 & 1 & 0\\
0 & 0 & 1 
\end{pmatrix}\,|\,z\,\in\mathbb{R}\right\}
\longrightarrow 
H_3 
\longrightarrow
\left\{\begin{pmatrix}
1 & x & 0\\
0 & 1 & y\\
0 & 0 & 1 
\end{pmatrix}\,|\,x,y\,\in\mathbb{R}\right\}
\longrightarrow
0,
$$
we know that 
the Heisenberg manifold is the total space of a principal
bundle with the structure group 
$$
U\left(1\right)\cong \mathbb{R}/\{2\ell\cdot k\,|\,k\in\mathbb{Z}\}
$$
and a torus
$$
T^2\cong \left(\mathfrak{g}_+\oplus\mathfrak{g}_-\right)/\Gamma_B
$$
as the base manifold.
We denote the projection map $H_3/\Gamma_{\ell}\to T^2$
by $\rho$.

\begin{prop}
The fibers of the bundle $\rho : H_3/\Gamma_{\ell}\to T^2$ are
totally geodesic and so 
the map 
$\rho^* : C^{\infty}\left(T^2\right)
\to C^{\infty}\left(H_3/\Gamma_{\ell}\right)$ 
commutes with the action of Laplacians on the each space.
\end{prop}

By this property the eigenvalues of the type (b) are those
coming from the base space $T^2$ through the map $\rho$.

\section{Zeta-regularized Determinant} 

Let ${\bf M}$ be an $n$-dimensional closed 
Riemannian manifold with the Laplacian
$\Delta_{\bf M}$, and denote the heat kernel
by $k_{\bf M}(t;x,y)$:
$$
\int_{\bf M} k_{\bf M}\left(t;x,y\right)
f\left(y\right)dy = e^{-t\Delta_{\bf M}}\left(f\right)(x).
$$ 
Then the Mellin transformation of the trace of the heat kernel, 
\begin{align*}
&\frac{1}{\Gamma\left(s\right)}\int_0^{\infty} 
\left(\int_{\bf M} k_{\bf M}(t;x,x)dx-1\right)t^{s-1}dt 
= \sum\limits_{\lambda\not=0}\frac{m_{\lambda}}{\lambda^{s}},\\
&m_{\lambda}\,\,\text{= multiplicity of the eigenvalue}\,\lambda,
\end{align*}
is meromorphically continued to the whole complex plane
with only poles of order one (at most)
at $s= n/2,\,n/2-1,\,\cdots$, especially
at $s=0$ it is holomorphic. 
We put this function as $\zeta_{\bf M}(s)$ and
call the spectral zeta-function of the Riemannian manifold ${\bf M}$. 
Then we can regard the value 
\[
e^{-\zeta_{\bf M}^{'}(0)}
\]
as a ``{\it determinant}'' (= product of non-zero eigenvalues) 
of the Laplacian $\Delta_{\bf M}$ acting on the space of functions 
orthogonal to the space of constant 
functions(\cite{RS}, \cite{QHS},\cite{Fo}) 
and call it as zeta-regularized determinant of the Laplacian. 
We denote it by $Det\,\, \Delta_{\bf M}$. Of course it can be
defined in a same way for more general elliptic operators.

In our case of the three dimensional Heisenberg manifold 
${\bf M}(\ell)=H_3/\Gamma_{\ell}$ we put
\allowdisplaybreaks{\begin{align*}
&\sum\limits_{\mu\in\Gamma_V(\ell)^{*},\mu\not=0}
\sum\limits_{m=0}^{\infty}
\left(Vol(\mathfrak{g}_+\oplus\mathfrak{g}_-/\Gamma_B)\cdot |\mu|\right)
\cdot e^{-t\{4\pi^2{\Vert\mu\Vert}^2+2\pi(2m+1)\Vert\mu\Vert\}}\\
&\qquad\qquad\qquad\qquad +\sum\limits_{\nu\in\Gamma_B^*}
e^{-4\pi^{2}t{\Vert\nu\Vert}^{2}}\\
&\qquad Z_{{\bf M}(\ell)}(t) = Z_V(t)+Z_{T^2}(t),
\end{align*}}
then
the second term is the trace of the heat kernel of the flat torus 
$T^2 \cong \left(\mathfrak{g}_+\oplus\mathfrak{g}_-\right)/\Gamma_B$, 
and so
\allowdisplaybreaks{\begin{align*}
\zeta_{{\bf M}(\ell)}(s)=&
\frac{1}{\Gamma(s)}\int_0^{\infty}(Z_{{\bf M}(\ell)}(t)-1)t^{s-1}dt\\
=\frac{1}{\Gamma(s)}\int_0^{\infty}&Z_V(t)t^{s-1}dt
+\frac{1}{\Gamma(s)}\int_0^{\infty}(Z_{T^2}(t)-1)t^{s-1}dt
= \zeta_V(s)+\zeta_{T^2}(s).
\end{align*}}   
Hence we have
\begin{prop}
\begin{equation*}
e^{-\zeta_{{\bf M}(\ell)^{'}(0)}}
=e^{-\zeta_V^{'}(0)}\cdot e^{-\zeta_{T^2}^{'}(0)}. 
\end{equation*}
\end{prop}

The value $\zeta_{T^2}^{'}(0)$ is given by the formula called
`` {\it Kroneker's second limit formula}'':
\begin{prop}{\bf (\cite{Be2}, \cite{Fo},\cite{QHS})}
$Det\,\,\Delta_{T^2} = e^{-\zeta_{T^2}^{'}(0)}= 
e^{-\frac{\pi}{3}}\,\left|\,
\prod\limits_{k=-\infty}^{\infty}
\left(1- e^{-2\pi\cdot |k|}\right)
\,\right|{^2}$.
\end{prop}
We give an elementary proof of this formula and expressions
of the zeta-reguralized  determinant 
for the three and four dimensional flat tori
in $\S 7$.

So in this section we only consider the value $\zeta_V^{'}(0)$.
The Mellin transform of the function $Z_V(t)$ is
\allowdisplaybreaks{\begin{align}\label{basic-formula}
&\frac{1}{\Gamma(s)}\int_0^{\infty}Z_V(t)t^{s-1}dt\\
&=4\ell\sum\limits_{n=1}^{\infty}\sum\limits_{m=0}^{\infty}
\frac{n}{(4\pi\ell)^{2s}(n^2+\frac{n}{4\pi\ell}(2m+1))^s}\notag\\
&=\frac{4\ell}{\Gamma(s)\Gamma(s-1)}
\cdot\frac{1}{(4\pi\ell)^{2s}}\cdot
\int_0^{\infty}\int_0^{\infty}
\frac{x+y}{e^{x+y}-1}\frac{x}{e^{\frac{x}{4\pi\ell}}-e^{-\frac{x}{4\pi\ell}}}
\frac{\quad (xy)^{s-2}}{x+y}dxdy.\notag
\end{align}}

Put 
$$f(x,y)
=\frac{x+y}{e^{x+y}-1}
\frac{x}{e^{\frac{x}{4\pi\ell}}-e^{-\frac{x}{4\pi\ell}}},
$$
then by the transformation
\[\left\{
\begin{array}{rl}
(x,y)\mapsto (u,v),& u=x\,\,\text{and}\,\,
v=\displaystyle{\frac{y}{x}}\,\, \text{on the domain}\,\, x>y\\
\quad &\\
(x,y)\mapsto (v,u),& u=y\,\,\text{and}\,\,
v=\displaystyle{\frac{x}{y}}\,\,\text{on the domain}\,\, x<y
\end{array}\right.
\]
we have
\allowdisplaybreaks{\begin{align*}
&\qquad\int_0^{\infty}\int_0^{\infty}f(x,y)\frac{\quad(xy)^{s-2}}{x+y}dxdy\\
&=\int_0^{1}\left(\int_0^{\infty}f\left(u,uv\right)u^{2s-4}du\right)
\frac{v^{s-2}}{1+v}dv
+\int_0^{1}\left(\int_0^{\infty}f\left(uv,u\right)u^{2s-4}du\right)
\frac{v^{s-2}}{1+v}dv\\
&=\int_0^{1}\left(\int_0^{\infty}
\left(f\left(u,uv\right)+f\left(uv,u\right)\right)u^{2s-4}du\right)
\frac{v^{s-2}}{1+v}dv\\
&=\int_0^{1}\left(\int_0^{\infty}
G\left(u,v\right)u^{2s-4}du\right)\frac{v^{s-2}}{1+v}dv,
\end{align*}}
where we put
$$
G(u,v)=f(u,uv)+f(uv,u).
$$


Let $g(x) = \displaystyle{\frac{x}{e^x-1}}$ and 
$h(x)=\displaystyle{\frac{2x}{e^x-e^{-x}}=\frac{x}{\sinh x}}$,\\
then
\begin{equation}\label{G(u,v)}
G(u,v) = 2\pi\ell\cdot
g\left(u\left(1+v\right)\right)
\left(h\left(\frac{u}{4\pi\ell}\right)+
h\left(\frac{uv}{4\pi\ell}\right)\right).
\end{equation}

Since the functions $g$ and $h$ are rapidly 
decreasing on the positive real axis, we have, 
\begin{prop}
For any integers $k$ and $l$ 
\begin{equation*}
\lim\limits_{u\to\infty}
\frac{\partial^{k+l}G}{\partial u^k\partial v^l}\left(u,v\right)=0
\end{equation*}
uniformly for $v\in [0,1]$.
\end{prop}

Next we consider the behavior of the functions 
$\displaystyle{\frac{\partial^{k+l}G}{\partial u^k\partial v^l}}$,
when $u\downarrow 0$.

Let 
\begin{equation}
g(x)=\frac{x}{e^x-1}
=\sum\limits_0^{\infty} \alpha_k x^k,\qquad |x|< 2\pi
\end{equation}
and
\begin{equation}
h(x)=\frac{2x}{e^x-e^{-x}}=\frac{x}{\sinh x}=
\frac{\sqrt{-1}x}{\sin \sqrt{-1}x}=\sum\limits_0^{\infty}
(-1)^k\beta_{2k} x^{2k},\qquad |x|<\pi.
\end{equation}
Note that $\alpha_0=1$, $\alpha_1=-1/2$, $\alpha_2=1/12$, 
$\alpha_{2i+1} = 0$ for $i=1,2,3,\cdots\cdots$  
and $\beta_0=1$, $\beta_2= 1/6$, $\beta_4=7/360$.
The coefficients are expressed in the following forms with
Bernoulli numbers $B_{2k}=\displaystyle{\frac{2(2k)!}{(2\pi)^{2k}}\zeta(2k)}$:
\begin{align*}
&\beta_{2k}=\frac{(2^{2k}-2)B_{2k}}{(2k)!}\\
&\alpha_{2k}= \frac{B_{2k}}{(2k)!}.
\end{align*} 


For $|u|<\pi$ and $v\in [0,1]$, the function $G(u,v)$
is expanded as follows:
\begin{equation*}
G(u,v)=2\pi\ell\cdot 
\sum\limits_{n=0}^{\infty}\left(
\sum\limits_{i+2j =n}(-1)^j\frac{\alpha_i\beta_{2j}}{(4\pi\ell)^{2j}}
(1+v)^i(1+v^{2j})\right)\cdot u^n 
=2\pi\ell\cdot\sum\limits_{n=0}^{\infty}P_n(v)u^n, 
\end{equation*}
where we denote the polynomial $P_n(v)$ 
$$P_n(v)=
\sum\limits_{i+2j =n}(-1)^j\frac{\alpha_i\beta_{2j}}{(4\pi\ell)^{2j}}
(1+v)^i(1+v^{2j}).
$$

\begin{prop} For $v \in[0,1]$
\begin{equation*}
\lim\limits_{u\downarrow 0}
\frac{\partial^{k+l}G}{\partial u^k\partial v^l}(u,v) 
= 2\pi\ell\cdot k!\cdot \frac{d^l P_{k}(v)}{dv^l}.
\end{equation*}
\end{prop}

Let $g(v,s)$ be a sufficiently many times differentiable function
defined on a domain in 
$\mathbb{R}\times \mathbb{C}$ 
including $[0,1]\times {\bf D}$, 
where 
${\bf D}=\left\{s\in\mathbb{C}\,|\, 
\mathfrak{Re}(s)> -\epsilon\right\}$ ($\epsilon > 0$ 
and fixed)
and $g$ is holomorphic 
on the domain ${\bf D}$ for each fixed
$t\in [0,1]$.

\begin{prop}\label{basic-res}
The function defined by the integral
\[
\int_0^1g(v,s)v^{s-2}dv
\]
has the Laurent expansion at $s=0$ as
\begin{equation}
\int_0^1g(v,s)v^{s-2}dv=
\frac{R_{-1}}{s}+R_0 +O(s),
\end{equation}
where $R_{-1}$ and $R_0$ are given by the formulas:

\allowdisplaybreaks{\begin{align*}
&R_{-1}=\frac{\partial g}{\partial v}(0,0)\\
&\text{and}\\
&R_0 = -\int_0^1\frac{\partial^2 g}{\partial v^2}(v,0)\log vdv
+\frac{\partial^2 g}{\partial s \partial v}(0,0)
+\frac{\partial g}{\partial v}(0,0)
-g(1,0).
\end{align*}}
\end{prop}

By applying this to the function of the form  
$\displaystyle{\frac{g(v,s)}{1+v}}$
we have the Laurent expansion of the function
\[
\int_0^1g(v,s)\frac{v^{s-2}}{1+v}dv
\]
as                    
\begin{cor}\label{residue2}
\allowdisplaybreaks{\begin{align*}
\int_0^1g(v,s)\frac{v^{s-2}}{1+v}dv&
=\left(\frac{\partial g}{\partial v}(0,0)-g(0,0)\right)
\frac{1}{s}\\
&-\int_0^1\frac{\partial^2 g}{\partial v^2}(v,0)\log vdv
+\int_0^1\frac{\partial g}{\partial v}(v,0)\log vdv
+\int_0^1\frac{g(v,0)}{1+v}dv\\
&+\frac{\partial^2 g}{\partial s \partial v}(0,0)
-\frac{\partial g}{\partial s}(0,0)
+\frac{\partial g}{\partial v}(0,0)-g(1,0) + 0(s).
\end{align*}} 
\end{cor}

When we restrict the variable $s$ in the domain
$\mathfrak{Re}(s) >3/2$, we have
\begin{equation}
\int_0^{\infty}G(u,v)u^{2s-4}du
=\frac{1}{(2s-3)(2s-2)(2s-1)2s}
\int_0^{\infty}\frac{\partial^4 G(u,v)}{\partial u^4}u^{2s}du.
\end{equation}

Now by Corollary (\ref{residue2}) when we put
\allowdisplaybreaks{\begin{align*}
&\int_0^{1}\left(\int_0^{\infty}
G(u,v)u^{2s-4}du\right)\frac{v^{s-2}}{1+v}dv\\
=&\frac{1}{(2s-3)(2s-2)(2s-1)2s}\int_0^{1}
\left(\int_0^{\infty}
\frac{\partial^4 G(u,v)}{\partial u^4}u^{2s}du\right)
\frac{v^{s-2}}{1+v}dv\\
&=\frac{1}{(2s-3)(2s-2)(2s-1)2s}\left\{ \frac{R_{-1}}{s}+R_0
  +O(1)\right\},
\end{align*}}
then
\begin{prop}  
\allowdisplaybreaks{\begin{align*}
&R_{-1}=0,\\
&R_0
=2\int_0^{\infty}
\frac{\partial^5 G(u,0)}{\partial v\partial u^4}\log udu
-2\int_0^{\infty}
\frac{\partial^4 G(u,0)}{\partial u^4}\log u du\\
&=-4\pi\ell\int_0^{\infty}
\frac{d^4}{du^4}
\left(h\left(\frac{u}{2}\right)^2\left(h\left(\frac{u}{4\pi\ell}\right)
+1 \right)\right)
\log u du.
\end{align*}}
\end{prop}
\begin{proof}
First we show
\allowdisplaybreaks{\begin{align*} 
R_{-1}&= \int_0^{\infty}\frac{\partial^5 G(u,0)}{\partial v\partial
  u^4}du
- \int_0^{\infty}\frac{\partial^4 G(u,0)}{\partial u^4}du\\
&=\frac{\partial^4 G(u,0)}{\partial v\partial u^3}
\left|\right._{0}^{\infty}
- \frac{\partial^3 G(u,0)}{\partial u^3}\left|\right._{0}^{\infty}\\
&=\frac{2\pi\ell}{(4\pi\ell)^2}\cdot g'(0)h''(0)
-\frac{2\pi\ell}{(4\pi\ell)^2}\cdot(1+v)\cdot g'(0)h''(0)|_{v=0}=0.
\end{align*}}
Next we calculate $R_0$:
\allowdisplaybreaks{\begin{align*}
R_0 =&-\int_0^1\left(\int_0^{\infty}\frac{\partial^6 G(u,v)}{\partial
  v^2\partial u^4}du\right) \log v dv
+\int_0^1\left(\int_0^{\infty}\frac{\partial^5 G(u,v)}{\partial
  v\partial u^4}du\right) \log vdv\notag\\
&+\int_0^1\left(\int_0^{\infty}
\frac{\partial^4 G(u,v)}{\partial u^4}du\right) \frac{1}{1+v}dv\notag\\
&+\int_0^{\infty}
\frac{\partial^5 G(u,0)}{\partial v\partial u^4}2\log udu
-\int_0^{\infty}\frac{\partial^4 G(u,0)}{\partial u^4}2\log u du\notag\\
&+\int_0^{\infty}\frac{\partial^5 G(u,0)}{\partial v\partial u^4}du
-\int_0^{\infty}\frac{\partial^4 G(u,1)}{\partial u^4}du\notag
\end{align*}}
\allowdisplaybreaks{\begin{align*}
&= \int_0^1\frac{\partial^5 G(0,v)}{\partial v^2\partial u^3}\log v dv
-\int_0^1\frac{\partial^4 G(0,v)}{\partial v\partial u^3}\log v dv
-\int_0^1\frac{\partial^3 G(0,v)}{\partial u^3}\frac{1}{1+v} dv\notag\\
&+2\int_0^{\infty}\frac{\partial^5 G(u,0)}{\partial v\partial u^4}\log udu
-2\int_0^{\infty}\frac{\partial^4 G(u,0)}{\partial u^4}\log u du\notag\\
&-\frac{\partial^4 G(0,0)}{\partial v\partial u^3}
+\frac{\partial^3 G(0,1)}{\partial u^3}\notag
\end{align*}}
\allowdisplaybreaks{
\begin{align*}
&=2\pi\ell\cdot 3!\cdot
\int_0^1\left(\frac{d^2 P_3(v)}{dv^2}-\frac{d P_3(v)}{dv}\right)
\log v - \frac{P_3(v)}{1+v}dv\notag\\
&+2\pi\ell\cdot 3!\cdot 
\left\{-\frac{d P_3(0)}{dv}+ P_3(1)\right\}\notag\\
&+2\int_0^{\infty}
\frac{\partial^5 G(u,0)}{\partial v\partial u^4} \log udu
-2\int_0^{\infty}
\frac{\partial^4 G(u,0)}{\partial u^4}\log u du\notag
\end{align*}
}
\allowdisplaybreaks{\begin{align*}
&=2\pi\ell\cdot 3!\cdot (-3)\cdot P_3(0)+2\pi\ell\cdot 
3!\cdot 3\cdot P_3(0)\notag\\
&+2\int_0^{\infty}
\frac{\partial^5 G(u,0)}{\partial v\partial u^4}\log udu
-2\int_0^{\infty}
\frac{\partial^4 G(u,0)}{\partial u^4}\log u du\notag
\end{align*}}
\allowdisplaybreaks{\begin{align*}
&=2\int_0^{\infty}
\frac{\partial^5 G(u,0)}{\partial v\partial u^4}\log udu
-2\int_0^{\infty}
\frac{\partial^4 G(u,0)}{\partial u^4}\log u du\notag\\
&=-4\pi\ell\int_0^{\infty}
\frac{d^4}{du^4}
\left(h\left(\frac{u}{2}\right)^2
\left(h\left(\frac{u}{4\pi\ell}\right)+1 \right)\right)
\log u du.\label{final form of G}
\end{align*}}
\end{proof}
Summing up we have
\allowdisplaybreaks{\begin{align*}
&\frac{1}{\Gamma(s)}\cdot\int_0^{\infty}Z_V(t)t^{s-1}dt
&=\frac{4\cdot\ell\cdot(s-1)}{\Gamma(s+1)^2}\cdot\frac{s^2}{(4\pi\ell)^{2s}}
\frac{1}{(2s-3)(2s-2)(2s-1)2s}\left\{R_0+ O(s)\right\}\\
&=\frac{\ell R_0}{3}s+O(s^2),
\end{align*}}
and the zeta-regularized determinant of the Laplacian
on the Heisenberg manifold $H_3/\Gamma_{\ell}$
is given by the formula:
\begin{thm}
\allowdisplaybreaks{\begin{align}
&Det\,\, \Delta_{H_3/\Gamma_{\ell}} = Det \,\,\Delta_{T^2}\cdot
e^{\frac{-\ell R_0}{3}}\notag\\
&=e^{-\frac{\pi}{3}}\,|\,
\prod\limits_{k=-\infty}^{\infty}(1- e^{-2\pi\cdot |k|})\,|{^2}
\cdot  e^{\frac{4\pi\ell^2}{3}\int_0^{\infty}
\frac{d^4}{du^4}\left(\left(\frac{{u}/{2}}{\sinh
    u/2}\right)^2\left(\frac{u/\left(4\pi\ell\right)}{\sinh
      \left(u/\left(4\pi\ell\right)\right)}+1\right)\right)\log udu}
\end{align}}
\end{thm}

\begin{cor}\label{limit of det 1}
When $\ell \to \infty$, $Det\,\, \Delta_{H_3/\Gamma_{\ell}}\to 0$. 
\end{cor}

\begin{proof}
Since 
\allowdisplaybreaks{\begin{align}
&\lim\limits_{\ell\to\infty}\int_0^{\infty}
\frac{d^4}{du^4}
\left\{
h\left(\frac{u}{2}\right)^2
\left(
h
\left(\frac{u}{4\pi\ell}
\right)
+1 \right)\right\}
\log u du\notag\\
&=2\int_0^{\infty}
\frac{d^4}{du^4}\left\{h\left(\frac{u}{2}\right)^2 \right\}
\log u du\label{limit form}
\end{align}}
we only determine the sign of this integral.
For this purpose we decompose the integral (\ref{limit form}) in the form 
\allowdisplaybreaks{\begin{align*}
&I=\int_0^{\infty}\frac{d^4}{du^4}
\left\{h\left(\frac{u}{2}\right)^2\right\}\log udu\\
&=\int_0^r\frac{d^4}{du^4}\left\{h\left(\frac{u}{2}\right)^2 - T_0
-T_1u-T_2u^2-T_3u^3\right\}\log udu\\
&+\int_r^{\infty}\frac{d^4}{du^4}
\left\{h\left(\frac{u}{2}\right)^2\right\}\log udu,
\end{align*}}
where
\allowdisplaybreaks{\begin{align*}
&T_0=h(0)^2=1,\\
&T_1=\frac{d}{du}h\left(\frac{u}{2}\right)^2_{|u=0}=0\\
&T_2=\frac{1}{2!}\frac{d^2}{du^2}
\left\{h\left(\frac{u}{2}\right)^2\right\}_{\left|u=0\right.}
=-\frac{1}{12}\\
&T_3=\frac{1}{3!}\frac{d^3}{du^3}
h\left(\frac{u}{2}\right)^2_{\left|\right.u=0}=0.
\end{align*}}
Then
\allowdisplaybreaks{\begin{align*}
&I=\frac{2}{r^{3}}
-\frac{1}{2r}-3!\left\{\int_0^r\left(h\left(\frac{u}{2}\right)^2 - 1
+\frac{1}{12}u^2\right)u^{-4}du
+\int_r^{\infty}
h\left(\frac{u}{2}\right)^2u^{-4}du\right\}.
\end{align*}}


Next we prove that $\displaystyle{h(u/2)^2-1+\frac{1}{12}u^2}$
takes positive values on the interval $[0,2]$.
For this purpose, recall that 
$h(u)= \displaystyle{\frac{u}{\sinh u}}$ is expanded for $u\in
(-\pi,\pi)$ as
\[
h(u)= \sum\limits_0^{\infty}(-1)^k\beta_{2k}u^{2k},
\]
with the coefficients
\[
\beta_{2k}= \frac{2(2^{2k}-2)}{(2\pi)^{2k}}\zeta(2k).
\]
Then 
\[\frac{\beta_{2k+2}}{\beta_{2k}}<
\frac{4-\frac{2}{2^{2n}}}{(2\pi)^2(1-\frac{2}{2^{2n}})}
\frac{\zeta(2k+2)}{\zeta(2k)}<\frac{3}{(2\pi)^2}<\frac{1}{12}.
\]
and
\[h(x)^2= \sum (-1)^k\gamma_{2k}x^{2k},
\]
with positive coefficients $\gamma_{2k}$ such that
\[
\gamma_{2k}= \sum\limits_{l=0}^{k}\beta_{2k-2l}\beta_{2l}.
\]
Hence
\[
\gamma_{2k}-\gamma_{2k+2}= 
\beta_0(\beta_{2n}-2\beta_{2n+2})
+\sum\limits_{l=1}^{k}\beta_{2l}(\beta_{2k-2l}-\beta_{2k+2-2l})
\]
so $\gamma_{2k}>\gamma_{2k+2}$ always.

{}From these estimates we see that the function
$\displaystyle{h(u/2)^2-1+\frac{1}{12}u^2}$
takes positive values on the interval $[0,2]$,
since $\displaystyle{h(u/2)^2-1+\frac{1}{12}u^2
=\sum\limits_{k=2}^{\infty}(-1)^k\gamma_{2k}
\left(\frac{u}{2}\right)^{2k}}$
is the form of the alternative sum consisting of
decreasing positive sequences when $u\in[0,2]$.  
Now we take $r=2$, then we see that $I<0$, hence we have the
desired result.
\end{proof}
\section{Heat kernel asymptotics}

As an application of the formula (\ref{basic-formula}) we calculate
the heat kernel asymptotics for the three dimensional 
Heisenberg manifolds, which are
given in terms of Bernoulli numbers.

We know that the heat kernel $k_{\bf M}(t;x,y)$ 
has the asymptotic expansion:
$$
k_{\bf M}(t;x,x)\sim \frac{1}{(4\pi t)^{n/2}}
\left\{c_0(x)+c_1(x) t +c_2(x) t^2+\cdots\cdots\right\},
\,\,t \downarrow 0, 
$$
$n$ = dimension of the manifold ${\bf M}$.
 
Of course the  coefficients  are given in terms of quantities coming from
metric tensors, but here we calculate the values
$$
{\bf c}_k=\int_{\bf M} c_k(x) dx
$$ 
explicitly by the method of analytic continuation by making use of 
the formula (\ref{basic-formula}), since we know that
the Mellin transform
of the heat kernel has poles of order one (at most)
at the points $n/2 - k,\, k = 0,1,2,\cdots\cdots$,
and the residue at the pole $n/2 - k$ is given by
the integral $\displaystyle{\int_{\bf M} c_k(x) dx}$. 

In our cases the trace of the heat kernel $Z_{H_3/\Gamma_{\ell}}(t) 
=\displaystyle{\int k_{H_3/\Gamma_{\ell}}(t;[g],[g])dg}$
is expanded as 
$$
\int_{H_3/\Gamma_{\ell}}k_{H_3/\Gamma_{\ell}}(t;[g],[g])dg 
= Z_V(t)+Z_{T^2}(t)
\sim \frac{1}{(4\pi t)^{3/2}}
({\bf c}_0+{\bf c}_1t+{\bf c}_2t^2+\cdots\cdots),
$$
and the second term $Z_{T^2}(t)$ is that corresponding to the two
dimensional flat torus. So that it is enough to consider the Mellin 
transform
of the first term $Z_V(t)$:
\allowdisplaybreaks{\begin{align}
&\frac{1}{\Gamma(s)}\int_0^{\infty}Z_V(t)t^{s-1}dt
=4\ell\sum\limits_{n=1}^{\infty}\sum\limits_{m=0}^{\infty}
\frac{n}{(4\pi\ell)^{2s}(n^2+\frac{n}{4\pi\ell}(2m+1))^s}\notag\\
&=\frac{4\ell}{\Gamma(s)\Gamma(s-1)}
\cdot\frac{1}{(4\pi\ell)^{2s}}\cdot
\int_0^{1}\left(\int_0^{\infty}
G(u,v)u^{2s-4}du\right)\frac{v^{s-2}}{1+v}dv\notag\\
&=\frac{4\ell}{\Gamma(s)\Gamma(s-1)}
\cdot\frac{1}{(4\pi\ell)^{2s}}\times\notag\\
&\times\left(\int_0^{1}\left(\int_0^{1}
G(u,v)u^{2s-4}du\right)\frac{v^{s-2}}{1+v}dv+\int_0^{1}\left(\int_1^{\infty}
G(u,v)u^{2s-4}du\right)\frac{v^{s-2}}{1+v}dv\right).
\label{basic-formula2}
\end{align}}
Here,   
\[G(u,v) = \frac{u(v+1)}{e^{u(v+1)}-1}\cdot
\left\{\frac{u}{e^{\frac{u}{4\pi\ell}}- e^{-\frac{u}{4\pi\ell}}}+ 
\frac{uv}{e^{\frac{uv}{4\pi\ell}}- e^{-\frac{uv}{4\pi\ell}}}\right\}. 
\]

Since the function defined by the integral
$$
\int_1^{\infty}G(u,v)u^{2s-4}du
$$
is holomorphic for any $s\in \mathbb{C}$ and so
the second term in the above expression (\ref{basic-formula2})
$$
\frac{4\ell}{\Gamma(s)\Gamma(s-1)}
\cdot\frac{1}{(4\pi\ell)^{2s}}
\int_0^{1}\left(\int_1^{\infty}G(u,v)u^{2s-4}du\right)
\frac{v^{s-2}}{1+v}dv
$$
has poles of order at most one at the points $s=1, 0,-1,-2, \cdots$.
Hence residues at these points must vanish.
Consequently it is enough to
consider the first term,
\begin{equation}\label{residue4}
\frac{4\ell}{\Gamma(s)\Gamma(s-1)}
\cdot\frac{1}{(4\pi\ell)^{2s}}\cdot
\int_0^{1}\left(\int_0^{1}
G(u,v)u^{2s-4}du\right)\frac{v^{s-2}}{1+v}dv,
\end{equation}
for calculating the residues at the poles $s= 3/2-k,
\,k=0,1,2,3\cdots$.

Then we have
\allowdisplaybreaks{\begin{align}
&\int_0^{1}\left(\int_0^{1}
G(u,v)u^{2s-4}du\right)\frac{v^{s-2}}{1+v}dv\notag\\
&=4\pi\ell\cdot
\int_0^{1}\left(\int_0^{1}\sum\limits_{i=0}^{\infty}\sum\limits_{j=0}^{\infty}
(-1)^j\alpha_i\beta_{2j}(u(1+v))^i
\left(
\left(\frac{u}{4\pi{\ell}}\right)^{2j}
+\left(\frac{uv}{4\pi{\ell}}\right)^{2j}
\right)
u^{2s-4}du\right)\frac{v^{s-2}}{1+v}dv\notag\\
&=4\pi\ell\cdot\sum\limits_{m=0}^{\infty}
\int_0^1\frac{1}{2s-3+m}
\sum\limits_{i+2j=m}(-1)^j\alpha_i\beta_{2j}\frac{1}{(4\pi \ell)^{2j}}
(1+v)^i(1+v^{2j})\frac{v^{s-2}}{1+v}dv.\label{residue3}
\end{align}}

To calculate the residues at the poles 
$s=3/2- k,\, k=0,1,2,\cdots\cdots$,
again it is enough to consider the terms corresponding to $m =
\text{even}= 2k$
in the sum (\ref{residue3}).

Let $W_{i,j}(s) = 
\displaystyle{\int_0^1(1+v)^{2i}(1+v^{2j})\frac{v^{s-2}}{1+v}dv}$,
then of course 
$W_{i,j}(s)$ is meromorphically continued to the whole complex plane
and we have
\begin{lem}\label{0-lem}
$W_{i,j}(s)$ is holomorphic at $s= 3/2-(i+j)$ and for $i>0$, 
$W_{i,j}(3/2-(i+j)) = 0$.
\end{lem}
\begin{proof}
Let $i>0$, then 
\begin{align*}
&\int_0^1(1+v)^{2i}(1+v^{2j})\frac{v^{s-2}}{1+v}dv\\
&=\sum\limits_{r=0}^{2i-1}\,_{2i-1}C_r
\left(\frac{1}{r+s-1}+\frac{1}{r+2j+s-1}\right).
\end{align*}
So at the point $s=3/2-(i+j)$ it takes the form
\allowdisplaybreaks{\begin{align*}
&\sum\limits_{r=0}^{2i-1}\,_{2i-1}C_r
\frac{1}{r+1/2 -(i+j)}\\
&+\sum\limits_{r=0}^{2i-1}
\,_{2i-1}C_{2i-1-r}\frac{1}{2i-1-r+2j+3/2-(i+j)-1}.
\end{align*}}
Hence we have $W_{i,j}(3/2-(i+j))=0$
\end{proof}

Now by Lemma (\ref{0-lem}), 
to calculate the residues of the function (\ref{residue2}) 
at the points $s= 3/2-k$, it is enough to consider the function
\begin{equation}
4\pi\ell \sum\limits_{k=0}^{\infty}\alpha_0\beta_{2k}
\frac{1}{(4\pi\ell)^{2k}}
\frac{1}{2s-3+2k}\int_0^1\frac{1+v^{2k}}{1+v}v^{s-2}dv.
\end{equation}

Let $W_k(s) 
=\displaystyle{\int_0^1\frac{1+v^{2k}}{1+v}v^{s-2}dv}$, then
the value at $s=3/2-k$, 
\begin{equation}\label{W_k}
W_k(3/2-k)={\bf W}_k 
\end{equation}
is given one by one by the following lemma
\begin{lem}
\allowdisplaybreaks{\begin{align*}
&W_k(3/2-k)\\
&=\int_0^1\frac{1+v^{2k}}{1+v}v^{s-2}dv_{|s=3/2-k}\\
&=\int_0^1\frac{v^{s-2}}{1+v}dv_{|s=3/2-k}+
\int_0^1\frac{v^{k-1/2}}{1+v}dv\\
&=\sum\limits_{r=0}^{k-1}\frac{(-1)^r}{r-k+1/2}
+(-1)^k\int_0^1\frac{1}{(1+v)\sqrt{v}}dv\\
&+\int_0^1\frac{v^k}{(1+v)\sqrt{v}}dv\\
&=\sum\limits_{r=0}^{k-1}\frac{(-1)^r}{r-k+1/2}
 +(-1)^k\frac{\pi}{2}
+2\sum_{r=0}^{k-1}\,_{k}C_r\left(\frac{k-r}{k(2r+1)}-J_r\right),
\end{align*}}
where we put 
$$
J_r=\int_0^1\frac{\theta^{2r}}{1+\theta^2}d\theta.
$$
$J_r$ is determined by the formula
\begin{equation*}
J_r = \sum\limits_{i=0}^{r-1}\,_{r}C_i
\left(\frac{r-i}{r(2i+1)}-J_i\right), J_0=\pi/4.
\end{equation*}
\end{lem}
Also this value ${\bf W}_k=W_k(3/2-k)$ is given by the formula:
\begin{equation*}
\sum\limits_{r=0}^{k-1}\frac{(-1)^r}{r-k+1/2}
 +(-1)^k\frac{\pi}{2} + 
\sum\limits_{j=0}^{\infty}\frac{(-1)^j}{j+k+1/2}.
\end{equation*}
Finally we have
\begin{prop}
The residue ${\bf c}_k =\displaystyle{\int_{H_3/\Gamma_{\ell}}c_k(g) dg}$
of the spectral zeta-function $\zeta_{H_3/\Gamma_{\ell}}(s)$ at 
the point $s=3/2-k,\,k -0,1,2,\cdots,$ is equal to
\begin{equation*}
\frac{4\pi}{\Gamma(3/2-k)\Gamma(1/2-k)}
\frac{\beta_{2k}}{2\cdot(4\pi\ell)^2}\cdot W_{k}(3/2-k).
\end{equation*}
\end{prop}
Note that 
$$
\beta_{2k}= \frac{2^{2k}-2}{(2k)!}B_{2k},
$$
where $B_{2k} = 
\displaystyle{\frac{2(2k)!}{(2\pi)^{2k}}\zeta(2k)}$ 
is the Bernoulli number.

\section{Five dimensional Heisenberg manifolds}

So far we only considered three dimensional cases and illustrated 
the procedure to calculate the determinants and heat kernel
asymptotics somehow precisely. These indicate that   
our method for calculating the determinants and heat kernel
asymptotics for higher dimensional Heisenberg manifolds
would also be valid. Even so the calculations and the
results are so complicated, here we 
state the results for the cases of 5-dimensional Heisenberg manifolds.

Let $\mathfrak{h}_{5}$ be the $5$-dimensional
Heisenberg Lie algebra:
\allowdisplaybreaks{\begin{align*}
\mathfrak{h}_5=&\mathfrak{g}_+\oplus\mathfrak{g}_-\oplus\mathfrak{z}
=\left\{ 
\begin{pmatrix}
0 & x_1 & x_2 & z \\
0 & 0   & 0   & y_1\\
0 & 0   & 0   & y_2\\
0 & 0   & 0   & 0 \\
\end{pmatrix}\,|\,x_i,\,y_i,\,z\in\mathbb{R}\right\}\\
&=\left\{ 
\begin{pmatrix}
0 & x_1 & x_2 & 0\\
0 & 0   & 0   & 0\\
0 & 0   & 0   & 0\\
0 & 0   & 0   & 0 \\
\end{pmatrix}\right\}\oplus
\left\{ 
\begin{pmatrix}
0 & 0 & 0 & 0 \\
0 & 0 & 0 & y_1\\
0 & 0 & 0 & y_2\\
0 & 0 & 0 & 0 \\
\end{pmatrix}\right\}\oplus
\left\{ 
\begin{pmatrix}
0 & 0 & 0 & z \\
0 & 0 & 0 & 0\\
0 & 0 & 0 & 0\\
0 & 0 & 0 & 0 \\
\end{pmatrix}\right\},
\end{align*}}
where $\mathfrak{z}$ 
is the center and $[\mathfrak{g}_+,\mathfrak{g}_-]=\mathfrak{z}$.

The corresponding Lie group $H_5$ is realized as
\begin{equation}
H_5 
=\left\{ 
\begin{pmatrix}
1 & x_1 & x_2 & z \\
0 & 1   & 0   & y_1\\
0 & 0   & 1   & y_2\\
0 & 0   & 0   & 1 \\
\end{pmatrix}\,|\,x_i,\,y_i,\,z\in\mathbb{R}\right\}.
\end{equation}

As in the three dimensional cases we only consider
a left invariant Riemannian metric defined from such an inner product
on the Lie algebra $\mathfrak{h}_5$ that
$X_1=\begin{pmatrix}
0 & 1 & 0 & 0\\
0 & 0 & 0 & 0\\
0 & 0 & 0 & 0\\
0 & 0 & 0 & 0\\
\end{pmatrix}$, 
$X_2=\begin{pmatrix}
0 & 0 & 1 & 0\\
0 & 0 & 0 & 0\\
0 & 0 & 0 & 0\\
0 & 0 & 0 & 0\\
\end{pmatrix}$, 
$Y_1=\begin{pmatrix}
0 & 0 & 0 & 0\\
0 & 0 & 0 & 1\\
0 & 0 & 0 & 0\\
0 & 0 & 0 & 0\\
\end{pmatrix}$, 
$Y_2=\begin{pmatrix}
0 & 0 & 0 & 0\\
0 & 0 & 0 & 0\\
0 & 0 & 0 & 1\\
0 & 0 & 0 & 0\\
\end{pmatrix}$ and
$Z_1=\begin{pmatrix}
0 & 0 & 0 & 1\\
0 & 0 & 0 & 0\\
0 & 0 & 0 & 0\\
0 & 0 & 0 & 0\\
\end{pmatrix}$\\
are orthonormal basis of $\mathfrak{h}_5$.

Also as in the three dimensional cases
let us take uniform discrete subgroups $\Gamma_{\ell}$ 
($\ell\in\mathbb{N}$) of the form
\begin{equation}
\Gamma_{\ell}=\left\{
\begin{pmatrix}
1 & m_1 & m_2 & \frac{k}{2\ell}\\
0 & 1   & 0   & n_1\\
0 & 0   & 1   & n_2\\
0 & 0   & 0   & 1 \\
\end{pmatrix}\left|\right.\, n_i,\,m_i ,\,k \in \mathbb{Z}\right\}.
\end{equation}

We identify $H_5$ and $\mathfrak{h}_5$ by the exponential map
\[
\exp:\mathfrak{h}_5 \to H_5
\]
\allowdisplaybreaks{
\begin{equation*}
\begin{array}{rl}
\exp:& g=\begin{pmatrix}
0 & x_1 & x_2 & z \\
0 & 0   & 0   & y_1\\
0 & 0   & 0   & y_2\\
0 & 0   & 0   & 0 \\
\end{pmatrix}\\
=&x_1X_1+x_2X_2+y_1Y_1+y_2Y_2+zZ_1 \mapsto 
\begin{pmatrix}
1 & x_1 & x_2 & z+1/2\sum x_iy_i\\
0 & 1   & 0   & y_1\\
0 & 0   & 1   & y_2\\
0 & 0   & 0   & 1 \\
\end{pmatrix}
\end{array}
\end{equation*}}
then the uniform discrete subgroup $\Gamma_{\ell}\subset H_5$ is
identified with a direct sum of two lattices $\Gamma_B\subset
\mathfrak{g}_+\oplus\mathfrak{g}_-$ 
and $\Gamma_V{(\ell)}\subset \mathfrak{z}$:
\allowdisplaybreaks{\begin{align*}
&\Gamma_B
=\left\{\begin{pmatrix}
0 & m_1 & m_2 & 0\\
0 & 0   & 0   & n_1\\
0 & 0   & 0   & n_2\\
0 & 0   & 0   & 0 \\
\end{pmatrix}\left|\right.\,m_i,\,n_i \in\mathbb{Z}\right\},\\
&\Gamma_V{(\ell)}
=\left\{\begin{pmatrix}
0 & 0 & 0 & \frac{k}{2\ell}\\
0 & 0   & 0   & 0\\
0 & 0   & 0   & 0\\
0 & 0   & 0   & 0 \\
\end{pmatrix}\left|\right.\,k\in\mathbb{Z}\right\},\\
&\exp\left(\Gamma_V(\ell)+\Gamma_B\right)=\Gamma_{\ell}
\end{align*}}

The heat kernel
on the whole group $H_5$ is given as the following form:
\allowdisplaybreaks{\begin{align}
&K_t(g,\tilde{g})=K_t(x,y,z;\tilde{x},\tilde{y},\tilde{z})\notag\\
&=(2\pi)^{-3}
\int_{-\infty}^{+\infty}e^{\sqrt{-1}\cdot\eta\cdot\left\{\tilde{z}-z+
\frac{1}{2}\left([\tilde{x},y]-[x,\tilde{y}]\right)\right\}}\cdot 
e^{-t\left|\eta\right|^2}\times\\
&\times\left(\frac{\left|\eta\right|}{2\sinh
    t\left|\eta\right|}\right)^2
\prod\limits_{i=1}^2
e^{-\frac{|\eta|}{4}\frac{\cosh t|\eta|}{\sinh
    t|\eta|}\cdot\left\{\left(x_i-\tilde{x}_i\right)^2+
\left(y_i-\tilde{y}_i\right)^2\right\}}d\eta.\notag
\end{align}}
where we regarded 
$\eta\in \mathfrak{z}^*=\mathbb{R}Z_1^*\cong\mathbb{R}$.
and $g=(x,y,z)=x_1X_1+x_2X_2+yY_1+y_2Y_2+zZ_1$, and similar to
$\tilde{g}$ through the exponential map.

Then the heat kernel $k_{H_5/\Gamma_{\ell}}(t;[g],[\tilde{g}])$ 
on the Heisenberg manifold 
${\bf M}_{\ell}= H_5/\Gamma_{\ell}$,
is expressed as :
\begin{equation}
k_{H_5/\Gamma_{\ell}}(t;[g],[\tilde{g}])=
\sum\limits_{\gamma\in \Gamma_{\ell}} K_t(\gamma\cdot g,\tilde{g}).
\end{equation}
and its trace is calculated in the following form:
\begin{thm}{\bf (\cite{Fu})}
\allowdisplaybreaks{\begin{align*}
&Z_{\ell}(t)\\
&=\int_{H_5/\Gamma_{\ell}}k_{H_5/\Gamma_{\ell}}(t;[g],[g])dg
=\sum\limits_{\gamma\in\Gamma_{\ell}}\int_{F_{\ell}}K_t(\gamma\cdot g,g)dg\\
&=\sum\limits_{\mu\in\Gamma_V(\ell)^{*},\mu\not=0}\sum\limits_{m_1=0}^{\infty}
\sum\limits_{m_2=0}^{\infty}
Vol(\mathfrak{g}_\oplus\mathfrak{g}_-/\Gamma_B)\Vert\mu\Vert^2
\cdot 
e^{-t\left\{\,4\pi^2{\Vert\mu\Vert}^2
+2\pi\sum\limits_{i=1}^2(2m_i+1)\Vert\mu\Vert\,\right\}}\\
&+\sum\limits_{\nu\in\Gamma_B^*}e^{-4\pi^{2}t{\Vert\nu\Vert}^{2}}\\
&=Z_V(t)+Z_{T^4}(t).
\end{align*}}
\end{thm}

Here $F_{\ell}$ denotes a fundamental domain of the uniform
discrete subgroup $\Gamma_{\ell}$ and 
$\Gamma_B^*$ and $\Gamma_V(\ell)^*$ 
are dual lattices, as before.

The second term $Z_{T^4}(t)$ is that corresponding to 
4-dimensional flat torus 
$(\mathfrak{g}_+\oplus\mathfrak{g}_-)/{\Gamma_B}$.

Let $\zeta_{{\bf M}_{\ell}}(s)$ be the function
\allowdisplaybreaks{\begin{align*}
\zeta_{{\bf M}_{\ell}}(s)=&
\frac{1}{\Gamma(s)}\int_0^{\infty}(Z(t)-1)t^{s-1}dt\\
=\frac{1}{\Gamma(s)}\int_0^{\infty}&Z_V(t)t^{s-1}dt
+\frac{1}{\Gamma(s)}\int_0^{\infty}(Z_{T^4}(t)-1)t^{s-1}dt
= \zeta_{V,\ell}(s)+\zeta_{T^4}(s),
\end{align*}}
then 
\begin{thm}
\[
Det\,\,\Delta_{H_5/\Gamma_{\ell}} = Det\,\,\Delta_{T^4}\cdot
e^{-\zeta_{V,\ell}^{\quad '}(0)},
\]
where $\displaystyle{{\zeta_{V,\ell}^{\quad '}(0)}}$  
is given in Proposition {\bf (}\ref{zeta_4-dim}{\bf )} below.
\end{thm}

Corresponding to the formula (\ref{basic-formula}) we have
an integral representation of the function $\zeta_{V,\ell}(s)$:
\allowdisplaybreaks{\begin{align}\label{basic-formula5}
&\zeta_{V,\ell}(s)=\frac{1}{\Gamma(s)}\int_0^{\infty}Z_V(t)t^{s-1}dt\\
&=8\ell^2\sum\limits_{n=1}^{\infty}\sum\limits_{m_1=0}^{\infty}
\sum\limits_{m_2=0}^{\infty}
\frac{n^2}{(4\pi\ell)^{2s}
\left(n^2+\frac{n}{4\pi\ell}\left(2m_1+1+2m_2+1\right)\right)^s}\notag\\
&=\frac{8\ell^2}{\Gamma(s)\Gamma(s-2)}
\cdot\frac{1}{(4\pi\ell)^{2s}}\cdot
\int_0^{\infty}\int_0^{\infty}
\frac{x+y}{e^{x+y}-1}
\left(\frac{x}{e^{\frac{x}{4\pi\ell}}-e^{-\frac{x}{4\pi\ell}}}\right)^2
\frac{\quad (xy)^{s-3}}{x+y}dxdy\notag\\
&=\frac{8\ell^2s^2(s-1)(s-2)} 
{\Gamma(s+1)^2\cdot (4\pi\ell)^{2s}}\cdot
\frac{4\pi^2\ell^2}{(2s-5)(2s-4)(2s-3)(2s-2)(2s-1)2s}\times\notag\\
&\times
\int_0^1\int_0^{\infty}\frac{d^6}{du^6}
\left\{g(u(1+v))\left(h\left(\frac{u}{4\pi\ell}\right)^2
+h\left(\frac{uv}{4\pi\ell}\right)^2\right)\right\}
u^{2s}\frac{v^{s-3}}{1+v}dudv.\notag
\end{align}}
\begin{prop}
\allowdisplaybreaks{\begin{align*}
\int_0^1\int_0^{\infty}\frac{d^6}{du^6}&
\left\{g(u(1+v))\left(h\left(\frac{u}{4\pi\ell}\right)^2
+h\left(\frac{uv}{4\pi\ell}\right)^2\right)\right\}
u^{2s}\frac{v^{s-3}}{1+v}dudv\\
=2\int_0^{\infty}\frac{d^6}{du^6}&
\left\{h\left(\frac{u}{2}\right)^3\left(\cosh \frac{u}{2}\right)
\left(h\left(\frac{u}{4\pi\ell}\right)^2+1\right)\right.\\
&\qquad\qquad\left.-\frac{1}{3}\left(\frac{u}{4\pi\ell}\right)^2
e^{-\frac{u}{2}}h\left(\frac{u}{2}\right)\right\}\log udu\\
=2\int_0^{\infty}\frac{d^6}{du^6}&
\left\{\left(\frac{\frac{u}{2}}{\sinh \frac{u}{2}}\right)^3 
\left(\cosh \frac{u}{2}\right)
\left(\left(
\frac{\frac{u}{4\pi\ell}}{\sinh \frac{u}{4\pi\ell}}\right)^2
+1\right)
\right.\\
&\qquad\qquad\left.-\frac{1}{3}
\left(\frac{u}{4\pi\ell}\right)^2
e^{-\frac{u}{2}}\frac{\frac{u}{2}}
{\sinh \frac{u}{2}}\right\} \log udu + 0(s).
\end{align*}}
\end{prop}
\begin{prop}\label{zeta_4-dim}
\allowdisplaybreaks{\begin{align*}
\zeta_{V,\ell}^{\,'}(0)& = -\frac{8}{15}\pi^2\ell^4\times\\
\times &\int_0^{\infty}\frac{d^6}{du^6}
\left\{\left(\frac{\frac{u}{2}}{\sinh \frac{u}{2}}\right)^3 \cosh u/2
\left(\left(
\frac{\frac{u}{4\pi\ell}}{\sinh \frac{u}{4\pi\ell}}\right)^2+1\right)
\right.\\
&\qquad\qquad\left.-\frac{1}{3}\left(\frac{u}{4\pi\ell}\right)^2
e^{-\frac{u}{2}}\frac{\frac{u}{2}}{\sinh \frac{u}{2}}\right\} \log udu.
\end{align*}}
\end{prop}
\begin{cor}
\allowdisplaybreaks{\begin{align*}
&\lim\limits_{\ell\to\infty}Det\,\, \Delta_{H_5/\Gamma_{\ell}}\\
&=\lim\limits_{\ell\to\infty}
Det\,\,\Delta_{T^4}\times e^{-\zeta_{V,\ell}^{'}(0)}=0.
\end{align*}}
\end{cor}
\begin{proof}
Put $I(\ell)$ 
\allowdisplaybreaks{\begin{align*}
&=\int_0^{\infty}\frac{d^6}{du^6}
\left\{\left(\frac{u/2}{\sinh u/2}\right)^3 \cosh u/2
\left(\left(\frac{u/(4\pi\ell)}{\sinh u/(4\pi\ell)}\right)^2+1\right)
\right.\\
&\qquad\qquad\left.-\frac{1}{3}\left(\frac{u}{4\pi\ell}\right)^2
e^{-u/2}\frac{u/2}{\sinh u/2}\right\} \log udu.
\end{align*}}
Since 
\begin{equation}\label{limit_integral}
\lim\limits_{\ell\to\infty}I(\ell)=2\int_0^{\infty}\frac{d^6}{du^6}
\left\{\left(\frac{u/2}{\sinh u/2}\right)^3 \cosh u/2 \right\}\log udu,
\end{equation}
it is enough to see the sign of this integral 
as we did before in Corollary (\ref{limit of det 1}) 
to determine the behavior of the
determinant $Det\,\, \Delta_{H_5/\Gamma_{\ell}}$ when $\ell\to\infty$.
Then, by a similar calculation in Corollary \ref{limit of det 1}
we have an expression of the integral 
\allowdisplaybreaks{\begin{align*}
&\int_0^{\infty}\frac{d^6}{du^6}
\left\{\left(\frac{u/2}{\sinh u/2}\right)^3 \cosh u/2 \right\}\log
udu\\
&=\frac{4!}{r^5}-\frac{1}{2r}-5!\int_0^{r}\left(
\left(\frac{u/2}{\sinh u/2}\right)^3 \cosh u/2 
-1+\frac{1}{2^4\cdot 3\cdot 5}u^4\right)u^{-6}du\\
&-5!\int_{r}^{\infty}\left(
\left(\frac{u/2}{\sinh u/2}\right)^3 \cosh u/2\right)u^{-6}du.
\end{align*}
}
So putting $r=2\cdot\sqrt[4]{3}$, we have
\allowdisplaybreaks{\begin{align*}
&\int_0^{\infty}\frac{d^6}{du^6}
\left\{\left(\frac{u/2}{\sinh u/2}\right)^3 \cosh u/2 \right\}\log
udu\\
&=-5!\int_0^{2{\sqrt[4]{3}}}\left(
\left(\frac{u/2}{\sinh u/2}\right)^3 \cosh u/2 
-1+\frac{1}{2^4\cdot 3\cdot 5}u^4\right)u^{-6}du\\
&-5!\int_{2{\sqrt[4]{3}}}^{\infty}\left(
\left(\frac{u/2}{\sinh u/2}\right)^3 \cosh u/2\right)u^{-6}du.
\end{align*}
}
It is clear that the second integrand 
takes always positive values on the
positive real axis. We can also prove that the integrand in the first
integral takes positive values on the positive real axis.
Since the coefficients of the Taylor expansion of the function
\begin{equation}
\cosh x -\left(1-\frac{1}{15}x^4\right)
\left(\frac{\sinh x}{x}\right)^3 =\sum\limits_{n=3}^{\infty}
a_{n}x^{2n}
\end{equation}
are given as 
\begin{equation} 
a_n=\frac{1}{(2n)!}+\frac{3^{2n-2}-1}{20\cdot(2n-1)!}
-\frac{3}{4}\frac{3^{2n+2}-1}{(2n+3)!}
\end{equation}
and all take positive values, which we can see from the 
expression of $a_n$ for $n\ge4$,
\allowdisplaybreaks{\begin{align*} 
a_n=\frac{1}{20\cdot(2n+3)!}&\left\{\left(3^{2n-2}-1\right)
\left((2n+3)(2n+2)(2n+1)(2n)-15\cdot 3^4\right)\right.\\
&\qquad+\left.20\cdot(2n+3)(2n+2)(2n+1)-15\cdot(3^4-1)\right\},\\
\qquad a_4= \frac{1}{20\cdot(11)!}&\left\{\left(3^{6}-1\right)
\left(11\cdot 10\cdot 9\cdot 8-15\cdot 3^4\right)\right.\\
&\qquad\quad+\left.20\cdot 11\cdot 10\cdot 9-15\cdot(3^4-1)\right\}
\end{align*}}
and for $n=3$, $a_3=\frac{4}{189}$.
{}From these facts we can prove the desired result.
\end{proof}

Finally we list the heat asymptotics
for five dimensional Heisenberg manifolds ${\bf M}_{\ell}$.
Let ${\bf c}_k$ be the coefficients of the asymptotic
expansion of the heat kernel $k_{H_5/\Gamma_{\ell}}(t;x,y)$ 
of the five dimensional Heisenberg manifolds
${\bf M}_{\ell}=H_5/\Gamma_{\ell}$:
\begin{equation*}
\int_{{\bf M}_{\ell}}k_{{\bf M}_{\ell}}(t;x,x)dx
\sim \frac{1}{(4\pi t)^{5/2}}\left\{
{\bf c}_0 +{\bf c}_1t+\cdots +{\bf c}_k t^k +\cdots\cdots\right\}.
\end{equation*}
Let us denote the Taylor expansion of the function 
$h(x)^2$ as
\begin{equation}
h(x)^2=\left(\frac{x}{\sinh x}\right)^2 =\sum\limits_{k=0}^{\infty}
(-1)^k\delta_k x^{2k},
\end{equation}
then $\delta_k$ is given by
\begin{equation}
\delta_k=\sum\limits_{j=0}^{k}\beta_{2j}\beta_{2(k-j)}
=\sum\limits_{j=0}^k \frac{4(2^{2(k-j)}-2)(2^{2j}-2)}{(2\pi)^{k}}
\zeta(2(k-j))\zeta(2j).
\end{equation}

\begin{prop}
\begin{equation}
{\bf c}_k =\frac{(2k-5)!!(2k-1)!!}{2^{2k+4}}
\frac{(-1)^k}{\pi^3\ell}\delta_{k}{\bf W}_k,
\end{equation}
where ${\bf W}_k$ is given in {\bf (}\ref{W_k}{\bf )}.
\end{prop}

 \section{A formula for product manifolds}

Let (${\bf M}, g$) and (${\bf N}, h$) 
be closed Riemannian manifolds, 
then the Laplacian $\Delta_{{\bf M}\times {\bf N}}$ 
on the product Riemannian manifold 
${\bf M} \times {\bf N}$
is of the form $\Delta_{\bf M}\otimes Id + Id \otimes\Delta_{\bf N}$ 
and the spectrum 
$Spec(\Delta_{{\bf M}\times {\bf N}})$ is given by
\begin{align*}
Spec(\Delta_{{\bf M}\times {\bf N}})=\{ \lambda_m +\mu_n\,|&\,
0=\lambda_0<\lambda_1\le\lambda_2\le \cdots \in Spec(\Delta_{\bf M}),\\
&\qquad\text{and}\,\,0=\mu_0<\mu_1\le \mu_2\le\cdots\in Spec(\Delta_{\bf N})\}.
\end{align*}

In this section we give a formula of the
zeta-regularized determinant of the Laplacian on the product
Riemannian manifold ${\bf M}\times {\bf N}$ 
in terms of the each value of the zeta-regularized
determinant and heat invariants.

The spectral zeta-function $\zeta_{{\bf M}\times {\bf N}}(s)$ 
for the product Riemannian manifold 
${\bf M}\times {\bf N}$ is given by
\[
\zeta_{{\bf M}\times {\bf N}}(s)
= \sum\limits_{m,n=0\,,\, (m,n)\not=0}^{\infty}
\frac{1}{(\lambda_m+\mu_n)^s}.
\]

We express this as
\allowdisplaybreaks{
\begin{align}
&\sum\limits_{m,n=0\,,\, (m,n)\not=0}^{\infty}
\frac{1}{(\lambda_m+\mu_n)^s}
=\frac{1}{\Gamma(s)}\sum\limits_{m=1}^{\infty}\frac{1}{\lambda_m^s}
\sum\limits_{n=0}^{\infty}\int_0^{\infty}
e^{-\left(1+\frac{\mu_n}{\lambda_m}\right)t}t^{s-1}dt
+\zeta_{\bf N}(s)\notag\\
&=\frac{1}{\Gamma(s)}
\sum\limits_{m=1}^{\infty}\frac{1}{\lambda_m^s}
\int_0^{\infty}
\left\{\sum\limits_{n=0}^{\infty} 
e^{-\frac{t}{\lambda_m}\cdot\mu_n}\right.\label{Q_(0)}\\
&\left.\qquad\qquad\qquad\qquad\qquad\qquad
-\left(\frac{\lambda_m}{4\pi t}\right)^{N/2}
\sum\limits_{i=0}^{[(N+M)/2]}
{\bf b}_i\cdot\left(\frac{t}{\lambda_m}\right)^i\right\}
e^{-t}t^{s-1}dt\notag\\
&\qquad\qquad+\frac{1}{\Gamma(s)}
\sum\limits_{m=1}^{\infty}\frac{1}{\lambda_m^s}
\int_0^{\infty}
\left(\frac{\lambda_m}{4\pi t}\right)^{N/2}
\sum\limits_{i=0}^{[(N+M)/2]}{\bf b}_i\cdot
\left(\frac{t}{\lambda_m}\right)^i
e^{-t}t^{s-1}dt +\zeta_{\bf N}(s)\notag\\
&=\mathcal{Q}_0(s) +\sum\limits_{i=0}^{[(N+M)/2]}
\frac{{\bf b}_i\cdot\Gamma(s+i-N/2)}{(4\pi)^{N/2}\cdot\Gamma(s)}
\cdot\zeta_{\bf M}(s+i-N/2)
+\zeta_{\bf N}(s)\label{basic decomposition},
\end{align}}
where $N= \dim {\bf N}$, $M = \dim {\bf M}$ and 
$\displaystyle{\{{\bf b}_i\}_{i=0}^{\infty}}$ 
are the coefficients of asymptotic
expansion of the trace
of the heat kernel $k_{{\bf N}}(t;x,y)$ on ${\bf N}$:
\begin{align}\label{Nasymptotic}
&k_{\bf N}(t)=\int_{{\bf N}}k_{{\bf N}}(t;x,x)dx\notag\\
&=\sum\limits_{n=0}^{\infty}e^{-t\mu_n}\sim
\left(\frac{1}{4\pi t}\right)^{N/2}
\left\{{\bf b}_0+ {\bf b}_1 t + {\bf b}_2 t^2 
+ {\bf b}_3 t^3+\cdots\cdots\right\}.
\end{align}

We also denote by $\{{\bf a}_i\}_{i=0}^{\infty}$ 
the coefficients of the asymptotic expansion of trace of
the heat kernel $k_{{\bf M}}(t;x,y)$
on the manifold ${\bf M}$: 
\begin{align}\label{Masymptotic}
&k_{\bf M}(t)=\int_{{\bf M}}k_{{\bf M}}(t;x,x)dx\notag\\
&=\sum\limits_{m=0}^{\infty} 
e^{-t\lambda_m}\sim \left(\frac{1}{4\pi t}\right)^{M/2}
\{{\bf a}_0+{\bf a}_1t+ {\bf a}_2t^2+ {\bf a}_3t^3+\cdots\cdots\}.
\end{align}

Since
\begin{equation*}
\left|\left\{
k_{\bf N}\left(\frac{t}{\lambda_m}\right)
-\left(\frac{\lambda_m}{4\pi t}\right)^{N/2}
\sum\limits_{i=0}^{[(N+M)/2]}
{\bf b}_i\cdot
\left(\frac{t}{\lambda_m}\right)^i\right\}\right|=
O\left(\left(\frac{t}{\lambda_m}\right)^{[(N+M)/2]+1-N/2}\right),
\end{equation*}
the first term $\mathcal{Q}_0(s)$ is holomorphic on the 
domain $\{s\in\mathbb{C}\,|\,
\mathfrak{Re}(s)> -1/2\}$, at least under this expression
for any case of the dimensions of the two manifolds. 
So we can put $s=0$ in (\ref{Q_(0)}) and we have 
\allowdisplaybreaks{
\begin{align}
&\mathcal{Q}_0(0)=0\notag\\
&\mathcal{Q}'_0(0)=
\sum\limits_{m=1}^{\infty}\int_0^{\infty}\left\{
k_{\bf N}\left(\frac{t}{\lambda_m}\right)
-\left(\frac{\lambda_m}{4\pi t}\right)^{N/2}
\sum\limits_{i=0}^{[(N+M)/2]}
{\bf b}_i\cdot
\left(\frac{t}{\lambda_m}\right)^i\right\}
e^{-t}t^{-1}dt\label{Q'_0(0)}.
\end{align}}

Next we put the second term as 
\allowdisplaybreaks{
\begin{align}
&\mathcal{Q}_1(s)=
\sum\limits_{i=0}^{[(N+M)/2]}
\frac{{\bf b}_i\cdot\Gamma(s+i-N/2)}{(4\pi)^{N/2}\cdot\Gamma(s)}
\cdot\zeta_{\bf M}(s+i-N/2)\label{Q_1all}
=\sum\limits_{i=0}^{[(N+M)/2]}{\mathfrak q}_i(s)\\
&=\sum\limits_{i=0}^{[(N+M)/2]}
\frac{{\bf b}_i}{(4\pi)^{N/2}\cdot\Gamma(s)}
\int_0^{\infty}
(k_{\bf M}(t)-1)t^{s+i-N/2-1}dt.\label{Q_1Mellin}
\end{align}}

\begin{rem}
By the asymptotic expansion {\bf (}\ref{Masymptotic}{\bf )} 
it is well known that the Mellin
transformation $\displaystyle{\int_0^{\infty}
(k_{\bf M}(t)-1)t^{s-1}dt}$ is meromorphically continued 
to the whole complex plane and has possible poles
of order one at points $M/2 - i,\, i=0,1,\cdots$
with the residue 
\[
\displaystyle{\frac{{\bf a}_i}{(4\pi)^{M/2}}}
\] when $\dim {\bf M}$ is odd, and when $\dim {\bf M}$ is even
then the residue at the pole $M/2 -i,\, i\not= M/2$ is
\[
\displaystyle{\frac{{\bf a}_i}{(4\pi)^{M/2}}},
\] 
and at $s=0$ the residue is
\[
\displaystyle{\frac{{\bf a}_{M/2}}{(4\pi)^{M/2}}}-1.
\]
\end{rem}

By the remark above 
we can describe the derivative at $s=0$ of the each term 
$${\mathfrak q}_i(s)
=\frac{{\bf b}_i\cdot\Gamma(s+i-N/2)}{(4\pi)^{N/2}\cdot\Gamma(s)}
\cdot\zeta_{\bf M}(s+i-N/2)
$$ 
in (\ref{Q_1all})
in terms of the spectral zeta-function $\zeta_{\bf M}(s)$. For this purpose
we consider four cases separately.

\begin{prop}\label{M*N e-o}
Let $\dim {\bf M}= M$ 
be even and  $\dim {\bf N}= N$  odd.
Since $\zeta_{\bf M}(s)$ is holomorphic at each point $i-N/2$ 
we have
\begin{align*}
&\mathcal{Q}_1'(0)= \sum\limits_{i=0}^{[(N+M)/2]}{\mathfrak q}_i'(0)\\
&=\sum\limits_{i=0}^{[(N+M)/2]} 
\frac{{\bf b}_i\cdot\Gamma(i-N/2)}{(4\pi)^{N/2}}\cdot\zeta_{\bf M}(i-N/2)
\end{align*}
\end{prop}

\begin{prop}\label{M*N o-e}
Let
$\dim {\bf M}= M $ 
be odd and  $\dim {\bf N}= N$  even.
Then
\allowdisplaybreaks{\begin{align*}
&for \,\,0\le i<N/2 \\
&\qquad\qquad {\mathfrak q}_i'(0)
=\frac{(-1)^{N/2-i}{\bf b}_i}{(4\pi)^{N/2}(N/2-i)!}\cdot
\zeta_{\bf M}'(i-N/2).\\
&for\,\,i=N/2\\
&\qquad\qquad {\mathfrak q}_{N/2}'(0)
=\frac{{\bf b}_{N/2}}{(4\pi)^{N/2}}\cdot\zeta_{\bf M}'(0).\\
&for \,\,N/2 <i \le [(N+M)/2]\\
&\qquad\qquad {\mathfrak q}_i'(0)
=\frac{{\bf b}_i(i-N/2-1)!}{(4\pi)^{N/2}}\cdot\zeta_{\bf M}(i-N/2).
\end{align*}}
Hence 
\allowdisplaybreaks{\begin{align*}
\mathcal{Q}_1'(0)&= \sum\limits_{i=0}^{N/2-1}
\frac{(-1)^{N/2-i}{\bf b}_i}{(4\pi)^{N/2}(N/2-i)!}\cdot
\zeta_{\bf M}'(0)
\\
&+\frac{{\bf b}_{N/2}}{(4\pi)^{N/2}}\cdot\zeta_{\bf M}'(0)\\
&+\sum\limits_{i=N/2+1}^{[(N+M)/2]}
\frac{{\bf b}_i}{(4\pi)^{N/2}(i-N/2-1)!}
\cdot\zeta_{\bf M}(i-N/2).
\end{align*}}
\end{prop}

\begin{prop}\label{M*N o-o}
Let both of $\dim {\bf M}= M $ 
and  $\dim {\bf N}= N$ be odd.
Then,
\allowdisplaybreaks{
\begin{align*}
&{\mathfrak q}_i'(0)=\frac{{\bf b}_i}{(4\pi)^{N/2}}
\left\{-\Gamma'(1)\cdot
\frac{{\bf a}_{\{(N+M)/2-i\}}}{(4\pi)^{M/2}}\right.\\
&\left.\qquad\qquad+
\left(\lim\limits_{s\to 0}\Gamma(i-N/2)\cdot
\zeta_{\bf M}(s+i-N/2)-
\frac{{\bf a}_{\{(N+M)/2-i\}}}{(4\pi)^{M/2}}\cdot\frac{1}{s}\right)
\right\}.
\end{align*}
}
Hence 
\allowdisplaybreaks{
\begin{align*}
&\mathcal{Q}_1'(0)=\sum\limits_{i=0}^{(N+M)/2}\frac{{\bf b}_i}{(4\pi)^{N/2}}
\left\{-\Gamma'(1)\cdot
\frac{{\bf a}_{\{(N+M)/2-i\}}}{(4\pi)^{M/2}}\right.\\
&\left.\qquad\qquad
+\left(\lim\limits_{s \to 0}
\Gamma(i-N/2)\cdot\zeta_{\bf M}(s+i-N/2)
-\frac{{\bf a}_{\{(N+M)/2-i\}}}{(4\pi)^{M/2}}\cdot\frac{1}{s}\right)
 \right\}
\end{align*}
}
\end{prop}   

\begin{prop}\label{M*N e-e} 
Let $\dim {\bf M}= M$ 
and $\dim {\bf N}= N$ be both even.
Then,
\allowdisplaybreaks{
\begin{align*}
&for \,\,0\le i< N/2\\
&\qquad\qquad {\mathfrak q}_i'(0)=\frac{{\bf b}_i}{(4\pi)^{N/2}}\left\{
-\Gamma'(1)\cdot\frac{{\bf a}_{\{(N+M)/2-i\}}}{(4\pi)^{M/2}}\right.\\
&\left.\qquad\qquad +
\lim\limits_{s \to 0}
\Gamma(s+i-N/2)\cdot\zeta_{\bf M}(s+i-N/2)
-\frac{{\bf a}_{\{(N+M)/2-i\}}}{(4\pi)^{M/2}}\cdot\frac{1}{s}
\right\},\\
& for\,\, i=N/2 \\
&\qquad\qquad q_{N/2}'(0)
=\frac{{\bf b}_{N/2}}{(4\pi)^{N/2}}\cdot\zeta_{\bf M}'(0)\,\, and\\
&for \,\, N/2< i\le (N+M)/2\\
&\qquad\qquad {\mathfrak q}_i'(0)=\frac{{\bf b}_i}{(4\pi)^{N/2}}
\left\{-\Gamma'(1)\cdot  
\frac{{\bf a}_{\{(N+M)/2-i\}}}{(4\pi)^{M/2}}\right.\\
&\left.\qquad\qquad
+\lim\limits_{s\to 0}(N/2-i-1)!\cdot\zeta_{\bf M}(s+i-N/2)
-\frac{{\bf a}_{\{(N+M)/2-i\}}}{(4\pi)^{M/2}}\cdot\frac{1}{s}
\right\}.
\end{align*}
} 
Hence,
\allowdisplaybreaks{
\begin{align*}
&\mathcal{Q}_1'(0)=\sum\limits_{i=0}^{N/2-1}
\frac{{\bf b}_i}{(4\pi)^{N/2}}\left\{
-\Gamma'(1)\cdot\frac{{\bf a}_{\{(N+M)/2-i\}}}{(4\pi)^{M/2}}\right.\\
&\left.\qquad\qquad +
\lim\limits_{s \to 0}\Gamma(s+i-N/2)\cdot\zeta_{\bf M}(s+i-N/2)
-\frac{{\bf a}_{\{(N+M)/2-i\}}}{(4\pi)^{M/2}}\cdot\frac{1}{s}
\right\}\\
&+\frac{{\bf b}_{N/2}}{(4\pi)^{N/2}}\cdot\zeta_{\bf M}'(0)\\
&+\sum\limits_{i=N/2+1}^{(N+M)/2}\frac{{\bf b}_i}{(4\pi)^{N/2}}
\left\{-\Gamma'(1)\cdot\frac{{\bf a}_{\{(N+M)/2-i\}}}{(4\pi)^{M/2}}\right.\\
&\left.\qquad\qquad
+\lim\limits_{s\to 0}(N/2-i-1)!\cdot\zeta_{\bf M}(s+i-N/2)
-\frac{{\bf a}_{\{(N+M)/2-i\}}}{(4\pi)^{M/2}}\cdot\frac{1}{s}
\right\}.
\end{align*}
}
\end{prop}

\begin{rem}
$-\Gamma'(1)=-\int_0^{\infty}e^{-t}\log t dt=
{\bf C_e}= 0.57721\cdots\cdots$ = Euler's constant.
\end{rem}

We can now write down an expression of the value
$Det\,\, \Delta_{{\bf M}\times {\bf N}}$ corresponding 
to the each case(\ref{M*N e-o}, \ref{M*N o-e}, \ref{M*N o-o} and
\ref{M*N e-e}) above.
Here we only
state for the case that both of $\dim {\bf M}$ and 
$\dim {\bf N}$ are even. In the next section
we state a special case of ${\bf N}= S^1$.
 
\begin{thm}\label{Det M*N e-e}
Let $\dim {\bf N}$ and $\dim {\bf M}$
be both even, then
\allowdisplaybreaks{\begin{align*}
&\qquad Det\,\, \Delta_{{\bf M}\times {\bf N}}
=e^{-\mathcal{Q}_0'(0)}\cdot e^{-\mathcal{Q}_1'(0)}\cdot Det\,\,\Delta_{\bf N}\\
&=Det\,\,\Delta_{\bf N}
\cdot
\prod\limits_{m=1}^{\infty}
e^{-\int_0^{\infty}\left\{
k_{\bf N}\left(\frac{t}{\lambda_m}\right)
-\left(\frac{\lambda_m}{4\pi t}\right)^{(N/2)}
\sum\limits_{i=0}^{(N+M)/2}
{\bf b}_i\cdot\left(\frac{t}{\lambda_m}\right)^i\right\}
e^{-t}t^{-1}dt
}\times  \\
&\times\prod\limits_{i=0}^{N/2-1}
e^
{
-\frac{{\bf C_e}\cdot{\bf b}_i\cdot{\bf a}_{\{(N+M)/2-i\}}}
{(4\pi)^{(N+M)/2}}
}
\cdot
\prod\limits_{i=0}^{N/2-1}
e^
{-\frac{{\bf b}_i}{(4\pi)^{(N/2)}}
\left\{\lim\limits_{s \to 0}\Gamma(s+i-N/2)
\cdot\zeta_{\bf M}(s+i-N/2)
-\frac{{\bf a}_{\{(N+M)/2-i\}}}{(4\pi)^{(M/2)}}\cdot\frac{1}{s}
\right\}
}\times\\
&\times 
(Det\,\,\Delta_{\bf M})^{\frac{{\bf b}_{N/2}}{(4\pi)^{(N/2)}}}
\cdot\prod\limits_{i=N/2+1}^{(N+M)/2}
e^
{-\frac{{\bf C_e}\cdot{\bf b}_{i}\cdot
{\bf a}_{\{(N+M)/2-i\}}}{(4\pi)^{(N+M)/2}}
}\times\\
&\times\prod\limits_{i=N/2+1}^{(N+M)/2}
e^
{-\frac{{\bf b}_i}{(4\pi)^{(N/2)}}
\left\{\lim\limits_{s\to 0}(N/2-i-1)!\cdot\zeta_{\bf M}(s+i-N/2)
-\frac{{\bf a}_{\{(N+M)/2-i\}}}{(4\pi)^{(M/2)}}\cdot\frac{1}{s}
\right\}
}.
\end{align*}
}
By interchanging ${\bf N}$ and ${\bf M}$ we have another
expression of $Det\,\, \Delta_{{\bf N}\times {\bf M}}$:
\allowdisplaybreaks{
\begin{align*}
&Det\,\, \Delta_{{\bf M}\times {\bf N}}
=Det\,\,\Delta_{\bf M}\cdot\prod\limits_{n=1}^{\infty}
e^{-\int_0^{\infty}\left\{
k_{\bf M}\left(\frac{t}{\mu_n}\right)
-\left(\frac{\mu_n}{4\pi t}\right)^{(M/2)}
\sum\limits_{i=0}^{(M+N)/2}
{\bf a}_i\cdot
\left(\frac{t}{\mu_n}\right)^i\right\}
e^{-t}t^{-1}dt}\times\\  
&\times\prod\limits_{i=0}^{M/2-1}
e^{-\frac{{\bf C_e}\cdot{\bf a}_i\cdot{\bf b}_{\{(N+M)/2-i\}}}{(4\pi)^{(M+N)/2}}}
\cdot\prod\limits_{i=0}^{M/2-1}
e^
{-\frac{{\bf a}_i}{(4\pi)^{(M/2)}}
\left\{\lim\limits_{s \to 0}\Gamma(s+i-M/2)\cdot\zeta_{\bf N}(s+i-M/2)
-\frac{{\bf b}_{\{(N+M)/2-i\}}}{(4\pi)^{(N/2)}}\cdot\frac{1}{s}
\right\}
}\times\\
&\times 
(Det\,\,\Delta_{\bf N})^{\frac{{\bf a}_{M/2}}{(4\pi)^{(M/2)}}}
\cdot\prod\limits_{i=M/2+1}^{(N+M)/2}
e^
{-\frac{{\bf C_e}\cdot{\bf a}_{i}\cdot{\bf b}_{\{(N+M)/2-i\}}}{(4\pi)^{(N+M)/2}}
}\times\\
&\times\prod\limits_{i=M/2+1}^{(N+M)/2}
e^
{-\frac{{\bf a}_i}{(4\pi)^{(M/2)}}
\left\{\lim\limits_{s\to 0}(M/2-i-1)!\cdot\zeta_{\bf N}(s+i-M/2)
-\frac{{\bf b}_{\{(N+M)/2-i\}}}{(4\pi)^{(N/2)}}\cdot\frac{1}{s}
\right\}
}.
\end{align*}
}
\end{thm}

\section{Product manifold ${\bf M}\times S^1$}

The formula we gave in the last section says 
that even in the product
manifold case the zeta-regularized determinant is
not expressed in a simple way in terms of each
zeta-regularized determinant. 
In the paper \cite{Fo} a formula 
of the zeta-regularized determinant
for a manifold of the type ${\bf M}\times S^1$ is given, 
in fact a formula
is derived for a higher order operator of the form
with the variable in $S^1$
being separated, by reducing the problem to a boundary value problem
on the manifold ${\bf M}\times [0,1]$
In this section we give a direct proof of this formula 
by restricting ourselves to the case of Laplacians on 
${\bf M}\times S^1$ by following the same line as we did in the last
section.  Here in this special case 
the Poisson summation formula allows us to
make an integration of the term 
(\ref{Q'_0(0)}) and we arrive
at a precise formula.

Now in this case the Laplacian $\Delta_{{\bf M}\times S^1}$ 
is given by 
$\displaystyle{\Delta_{{\bf M}\times S^1}
=\Delta_{\bf M} - \frac{d^2}{dx^2}}$, 
where we regard
$S^1\cong\mathbb{R}/(2\pi\ell\cdot\mathbb{Z})$, 
$x\in\mathbb{R},\,\ell>0$ and
the spectrum
are 
\[
\left\{ \lambda_n+ \left(\frac{k}{\ell}\right)^2\,\left|\, 
0=\lambda_0<\lambda_1\le\lambda_2,\,\cdots,\,
\text{ are the spectrum of}\, {\bf M},
\,\text{and}\,k
\in\mathbb{Z}\right.\right\}.
\]

\begin{thm}
\begin{equation}\label{general_product_formula}
Det\,\,\Delta_{{\bf M}\times S^1}
=4\pi^2\ell^2 
{\bf C}_{\bf M}
\cdot \prod\limits_{n=1}^{\infty}
\left|\left(1-e^{-2\pi\ell\sqrt{\lambda_n}}\right)\right|^2,
\end{equation}
where the constant ${\bf C}_{\bf M}$ is given by
\end{thm}
\allowdisplaybreaks{
\begin{enumerate}
\item {\it when}
$\dim {\bf M} = \,{\it even}\, = 2m$, {\it then} 
\[{\bf C}_{\bf M} = e^{2\pi\ell\cdot\zeta_{\bf M}(-1/2)}
\]
\item {\it when} 
$\dim {\bf M} = \,{\it odd}\, = 2m+1$, {\it then} 
\[\log {\bf C}_{\bf M} = -\sqrt{\pi}\ell
\left\{
\lim\limits_{s\to 0}\left(
2\sqrt{\pi}\cdot\zeta_{\bf M}(s-1/2)
+\frac{{\bf a}_{(M+1)/2}}{(4\pi)^{(M+1)/2}}\frac{1}{s}
\right)
+\frac{\Gamma'(1)}{(4\pi)^{M/2}}{\bf a}_{(M+1)/2}
\right\}.
\]
\end{enumerate}}
Here ${\bf a}_{(M+1)/2}$ denotes $(M+1)/2$-th 
coefficient of the asymptotic expansion
of the heat kernel $Z_{\bf M}(t)$ on ${\bf M}$.
\begin{proof}
We express the 
spectral zeta-function $\zeta_{{\bf M}\times S^1}(s)$ as
\allowdisplaybreaks{
\begin{align*}
\zeta_{{\bf M}\times S^1}(s)&=\sum\limits_{n=1}^{\infty}
\sum\limits_{k\in\mathbb{Z}}
\frac{1}{\left(\lambda_n+\left(\frac{k}{\ell}\right)^2\right)^s}
+2\ell^{2s}\cdot\zeta(2s)\\
&=\frac{1}{\Gamma(s)}
\sum\limits_{n=1}^{\infty}
\frac{1}{\lambda_n^s}\int_0^{\infty}
\sum\limits_{k\in\mathbb{Z}}
e^{-\left(1+\left(\frac{k}{\ell}\right)^2
\frac{1}{\lambda_n}\right)x}x^{s-1}dx
+2\ell^{2s}\cdot\zeta(2s).
\end{align*}  
}
Then by using the Poisson's summation formula this equals to
\allowdisplaybreaks{
\begin{align*}
&\frac{1}{\Gamma(s)}
\sum\limits_{n=1}^{\infty}
\frac{1}{\lambda_n^s}\int_0^{\infty}
\sum\limits_{k\in\mathbb{Z}}
\sqrt{\frac{\pi \lambda_n \ell^2}{x}}
e^{-\frac{\pi^2 \lambda_n k^2\ell^2}{x}}
e^{-x}x^{s-1}dx + 2\ell^{2s}\cdot\zeta(2s)\\
&=\frac{1}{\Gamma(s)}
\sum\limits_{n=1}^{\infty}
\frac{1}{\lambda_n^s}\int_0^{\infty}
\sum\limits_{k\in\mathbb{Z},\,k\not=0}
\sqrt{\frac{\pi \lambda_n \ell^2}{x}}
e^{-\frac{\pi^2 \lambda_n k^2 \ell^2}{x}}e^{-x}x^{s-1}dx\\
&+\frac{1}{\Gamma(s)}
\sum\limits_{n=1}^{\infty}
\frac{1}{\lambda_n^s}\int_0^{\infty}
\sqrt{\frac{\pi \lambda_n \ell^2}{x}}e^{-x}x^{s-1}dx
+2\ell^{2s}\cdot\zeta(2s)\\
&=\mathcal{P}_0(s)
+\frac{\sqrt{\pi}\ell\cdot\Gamma(s-1/2)}{\Gamma(s)}
\sum\limits_{n=1}^{\infty}\frac{1}{\lambda_n^{s-1/2}}
+2\ell^{2s}\cdot\zeta(2s),
\end{align*}
}
where we put
\allowdisplaybreaks{
\begin{align*}
&\mathcal{P}_0(s)=\frac{1}{\Gamma(s)}
\sum\limits_{n=1}^{\infty}
\frac{1}{\lambda_n^s}\int_0^{\infty}
\sum\limits_{k\in\mathbb{Z},\,k\not=0}
\sqrt{\frac{\pi \lambda_n \ell^2}{x}}
e^{-\frac{\pi^2 \lambda_n k^2 \ell^2}{x}}e^{-x}x^{s-1}dx\\
&=\frac{2\sqrt{\pi}\ell}{\Gamma(s)}
\sum\limits_{n=1}^{\infty}
\frac{1}{\lambda_n^{s-1/2}}\int_0^{\infty}
\sum\limits_{k=1}^{\infty}
e^{-\frac{\pi^2 \lambda_n k^2 \ell^2}{x}}e^{-x}x^{s-3/2}dx.
\end{align*}
}
We denote the second term 
$\displaystyle{\frac{\sqrt{\pi}\ell\cdot\Gamma(s-1/2)}{\Gamma(s)}
\sum\limits_{n=1}^{\infty}\frac{1}{\lambda_n^{s-1/2}}}
=\frac{\sqrt{\pi}\ell\cdot\Gamma(s-1/2)}{\Gamma(s)}\cdot
\zeta_{\bf M}(s-1/2)$
by $\mathcal{P}_1(s)$.

Since the sum $\displaystyle{\sum\limits_{k=1}^{\infty}
e^{-\frac{\pi^2 \lambda_n k^2 \ell^2}{x}}}$ in the integrand 
satisfies the asymptotics
\[
\sum\limits_{k=1}^{\infty}
e^{-\frac{\pi^2 \lambda_n k^2\ell^2}{x}}
=O\left( \left( \frac{x}{\lambda_n} \right)^{N}\right)
\]
for any $N\in\mathbb{N}$, or the coefficients ${\bf b}_i$ are all zero
 except ${\bf b}_0= 2\pi\ell$,
 the function $\mathcal{P}_0(s)$ is holomorphic on
the whole complex plane and we have
\allowdisplaybreaks{\begin{align*}
&\mathcal{P}_0(0)=0,\\
&\mathcal{P}_0'(0)=2\sqrt{\pi}\ell
\sum\limits_{n=1}^{\infty}\sqrt{\lambda_n}
\int_0^{\infty}
\sum\limits_{k=1}^{\infty}
e^{-\frac{\pi^2 \lambda_n k^2 \ell^2}{x}}e^{-x}x^{-3/2}dx.
\end{align*}}

To calculate the value $\mathcal{P}_0'(0)$ 
recall a formula of 
a {\it modified Bessel function} 
$K_{1/2}(z)=K_{-1/2}(z)$ (see \cite{AAR}):
\begin{equation}\label{Bessel_1/2}
\int_0^{\infty}
e^{-\left(t+\frac{z^2}{4t}\right)}\frac{1}{\sqrt{t}}dt 
= \sqrt{\pi}e^{-z}= \sqrt{2z}K_{1/2}(z). 
\end{equation}
Now we have
\[
\mathcal{P}_0'(0)=-2\sum\limits_{n=1}^{\infty}
\log\left(1-e^{-2\pi\ell\sqrt{\lambda_n}}\right).
\]
Hence from Propositions \ref{M*N e-o} and \ref{M*N o-o}
the determinant $Det\,\,\Delta_{{\bf M}\times S^1}$
is of the form
\allowdisplaybreaks{
\begin{equation}\label{prod_formula}
Det\,\,\Delta_{{\bf M}\times S^1}
=4\pi^2\ell^2 
{\bf C}_{\bf M}
\cdot \prod\limits_{n=1}^{\infty}
\left|\left(1-e^{-2\pi\ell\sqrt{\lambda_n}}\right)\right|^2,
\end{equation}
}
where the constant ${\bf C}_{\bf M}$ is given by
\allowdisplaybreaks{
\begin{enumerate}
\item when
$\dim {\bf M}$ is even, then 
\[{\bf C}_{\bf M} = e^{2\pi\ell\cdot\zeta_{\bf M}(-1/2)}
\]
\item when 
$\dim {\bf M} = M$ is odd, then 
\[\log {\bf C}_{\bf M} = -\sqrt{\pi}\ell
\left\{
\lim\limits_{s\to 0}\left(2\sqrt{\pi} 
\cdot\zeta_{\bf M}(s-1/2)
+\frac{{\bf a}_{(M+1)/2}}{(4\pi)^{(M+1)/2}}\frac{1}{s}\right)
+\frac{\Gamma'(1)}{(4\pi)^{M/2}}{\bf a}_{(M+1)/2}\right\}.
\]
\end{enumerate}}
\end{proof}

\begin{example}
As an application of our formula 
{\bf (}\ref{general_product_formula}{\bf )} 
we give an expression of $Det\,\, \Delta_{S^2\times S^1}$.

For the standard 2-dimensional sphere $S^2$ the spectral 
zeta-function is 
$$
\zeta_{S^2}(s)=\sum\limits_{k=1}^{\infty}\frac{2k+1}{k^s(k+1)^s},
$$
which converges for $\frak{Re}(s)>1$.
We rewrite this as
\begin{align}
&\zeta_{S^2}(s)=\frac{1}{2^s}+
\sum\limits_{k=2}^{\infty}\frac{1}{k^{2s-1}}
\left\{\left(1+\frac{1}{k}\right)^{-s}
+\left(1-\frac{1}{k}\right)^{-s}\right\}\notag\\
&=\frac{1}{2^s}+
\sum\limits_{k=2}^{\infty}
\frac{1}{k^{2s-1}}\sum\limits_{m=0}^{\infty}2d_{2m}(-s)
\left(\frac{1}{k}\right)^{2m}\notag\\
&=\frac{1}{2^s}+
2\sum\limits_{m=0}^{2m\leq n}d_{2m}(-s)(\zeta(2s-1+2m)-1)\label{S2zeta}\\
&\qquad\qquad\qquad 
+2\sum\limits_{2m>n}^{\infty}\sum\limits_{k=2}^{\infty}
d_{2m}(-s)\frac{1}{k^{2m-n}}\cdot\frac{1}{k^{2s-1+n}},\notag
\end{align}
where we used the expansion 
$\displaystyle{(1+z)^{\alpha}=\sum d_m(\alpha)z^m}$, 
for $|z|<1$. Note that for $\alpha > 0$
this series converges for $-1\leq z \leq 1$.
Then for $\frak{Re}(s) >(2-n)/2$, by the estimate 
\begin{align*}
&\sum\limits_{2m>n}^{\infty}\sum\limits_{k=2}^{\infty}
|d_{2m}(-s)\frac{1}{k^{2m-n}}\frac{1}{k^{2s-1+n}}|
\leq \sum\limits_{2m>n}^{\infty}|d_{2m}(-s)|\frac{1}{2^{2m-n}}\cdot
\sum\limits_{k=2}^{\infty}\frac{1}{k^{2\frak{Re}(s)-1+n}}
\end{align*}
and the functional relation for the Riemann $\zeta$-function
the expression (\ref{S2zeta}) gives us the analytic continuation
of $\zeta_{S^2}(s)$ to the complex plane of $\frak{Re}(s)> (2-n)/2$
for each $n>0$. So we can put $s=-1/2$ in (\ref{S2zeta}) and
we have an expression of $\zeta_{S^2}(-1/2)$:
\begin{align*}
&\zeta_{S^2}(-1/2)\\
&=\sqrt{2}+2\left\{\zeta(-2)-1\right\}+
2d_2(1/2)\left\{\zeta(0)-1\right\}\\
&+2\sum\limits_{m=3}^{\infty}d_{2m}(1/2)(\zeta(2m-2)-1)\\
&=\sqrt{2}
-2\sum\limits_{m=0}^{\infty}d_{2m}(1/2)
+2\sum\limits_{m=0}^{\infty}d_{2m}(1/2)(\zeta(2m-2))\\
&=-\sum\limits_{m=0}^{\infty}\frac{(4m)!}{2^{4m-1}(4m-1)((2m)!)^2}\zeta(2m-2).
\end{align*}
This is also expressed as
\begin{equation}
\zeta_{S^2}(-1/2)=\frac{4}{9\pi}\int_0^{\infty}\int_0^{\infty}
\frac{\partial^2}{\partial x^2}\left(
\frac{\partial^2}{\partial y^2}\left(\frac{x+y}{e^{x+y}-1}\right) 
e^{-x}\right)\cdot\frac{1}{\sqrt{xy}}dxdy.
\end{equation}
So, finally we have
\begin{equation}
Det\,\,\Delta_{S^2\times S^1}=
4\pi^2\ell^2 
\prod\limits_{m=1}^{\infty}
e^{-\pi\ell\frac{(4m)!}{2^{4m-2}(4m-1)((2m)!)^2}\zeta(2m-2)}
\cdot \prod\limits_{n=1}^{\infty}
\left|\left(1-e^{-2\pi\ell\sqrt{k(k+1)}}\right)^{2k+1}\right|^2.
\end{equation}
\end{example}
\smallskip

Again by applying Propositions \ref{M*N o-o} and \ref{M*N o-e} 
we have an alternative representation
of the determinant $Det\,\, \Delta_{{\bf M}\times S^1}$.
\begin{cor}
When $\dim {\bf M}$ is odd, then
\allowdisplaybreaks{
\begin{align}
&Det\,\, \Delta_{{\bf M}\times S^1}\\
&=\prod\limits_{k=1}^{\infty}
e^{-2\int_0^{\infty}
\left\{
k_{\bf M}\left(\frac{\ell^2}{k^2}x\right)
-\left(\frac{k^2}{4\pi\ell^2 x}\right)^{M/2}
\sum\limits_{j=0}^{(M+1)/2}
{\bf a}_j\cdot\left(\frac{\ell^2}{k^2}x\right)^j
\right\}
e^{-t}t^{-1}dt}\times\notag\\
&\times\prod\limits_{i=0}^{[M/2]}
e^{-\frac{2}{(4\pi\ell^2)^{M/2}}\cdot
{\bf a}_i\cdot\ell^{2i}\cdot\Gamma(i-M/2)\cdot\zeta(2i-M)}\times\notag\\
&\times
e^{-\frac{2}{(4\pi\ell^2)^{M/2}}{\bf a}_{[(M+1)/2]}\cdot\ell^{M+1}
\cdot\left\{
\sqrt{\pi}
\left(\log \ell -{\bf C_e}/2 +1/2\Gamma'(1/2)\right)
\right\}}
\times\notag\\
&\times Det\,\,\Delta_{\bf M}.\notag
\end{align}
}
Here we used the formula $\zeta(s)=1/(s-1)+ {\bf C_e}s + O((s-1)^2)$.

When $\dim {\bf M}$ is even, then
\allowdisplaybreaks{
\begin{align}
&Det\,\, \Delta_{{\bf M}\times S^1}\\
&=\prod\limits_{k=1}^{\infty}
e^{-2\int_0^{\infty}
\left\{
k_{\bf M}(\frac{\ell^2}{k^2}x)
-\left(\frac{k^2}{4\pi\ell^2 x}\right)^{M/2}
\sum\limits_{j=0}^{(M+1)/2}
{\bf a}_j\cdot\left(\frac{\ell^2}{k^2}x\right)^j
\right\}
e^{-t}t^{-1}dt}\times\notag\\
&\times\prod\limits_{i=0}^{M/2-1}
e^
{-\frac{4{\bf a}_i\cdot\ell^{2i}}{(4\pi\ell^2)^{M/2}}  
\frac{(-1)^{M/2-i}}{(M/2-i)!}
\cdot\zeta'(2i-M)}
\times\notag\\
&\times 2\pi\ell\cdot Det\,\,\Delta_{\bf M}.\notag
\end{align}
}
Note that $\zeta(-2k)=0$ for $k=1,2,\,\cdots$.
\end{cor}

\begin{rem}
Our formula {\bf (}\ref{general_product_formula}{\bf )} 
is of course a special case 
of formulas given in the paper \cite{Fo} for more general
elliptic operators on the product manifolds ${\bf M}\times S^1$.
However here we gave an expression of the constant ${\bf C}_{\bf M}$,
although the formula itself is not a computable form, especially
for ${\bf M}$ being odd dimensional. To obtain a further information
we must specify the manifolds ${\bf M}$. So in the next section 
we give a more precise form of this factor ${\bf C}_{\bf M}$
for some flat tori. 
\end{rem}

\section{Flat tori}

In the last two sections we considered the zeta-regularized
determinant for manifolds of a product
form as a Riemannian manifold.  
In this section we deal with the case that the manifolds are two, three
and four dimensional flat tori, which are not always
of a product form of lower dimensional tori as Riemannian
manifolds.

We know by a similar calculation as we showed in $\S 2$ and $\S 4$
that the zeta-regularized 
determinant of ($2n+1$)-dimensional Heisenberg manifolds 
are always of the product
form with a factor 
which is the zeta-regularized determinant 
of a $2n$-dimensional torus.  So of course 
it is required to determine the zeta-regularized
determinant of flat tori to complete the calculation for Heisenberg manifolds. 
In this section we give an expression of it for two, 
three and four dimensional flat tori. Although our expressions 
are not of a computable form within a finite step, the expression for two
dimensional cases are given by the famous limit formula of Kroneker as
we cited in $\S 2$, and higher cases correspond to a generalization of 
this limit formula. The structure of a generalization was already
stated and discussed focusing in their functional relations 
in the papers \cite{Be1} and \cite{Be2} for more general
Dirichlet series than Epstein zeta-functions which are of our cases.
Here we treat with the typical Epstein zeta-functions of 
two, three and four variables.
Since it is enough for our purpose to give an explicit analytic
continuation of the functions 
from a left half region in the complex plane 
to a region including zero, we give them based on the Jacobi identity
and the Mellin transformation in a quite elementary way.
For this purpose we fix the 
flat tori in the following way.

Let ${\bf e}_1,\,{\bf e}_2,\,{\bf e}_3,\,{\bf e}_4$ 
be the standard orthonormal basis on $\mathbb{R}^4$ and 
we fix a basis $\{{\bf u}_1,{\bf u}_2,{\bf u}_3,{\bf u}_4\}$ of the
following form 
\allowdisplaybreaks{
\begin{align*}
&{\bf u}_1 ={\bf e}_1,\, {\bf u}_2
=a_{1,2}{\bf e}_1+a_{2,2}{\bf e}_2,\,
\text{we put this}\,=A{\bf e}_1+B{\bf e}_2,\\
&{\bf u}_3=a_{1,3}{\bf e}_1 + a_{2,3}{\bf e}_2 + a_{3,3}{\bf e}_3 \\
&{\bf u}_4=a_{1,4}{\bf e}_1 + a_{2,4}{\bf e}_2 + a_{3,4}{\bf e}_3
+a_{4,4}{\bf e}_4\,(a_{2,2},\,a_{3,3},\,a_{4,4}\,>\,0).
\end{align*}
}
\allowdisplaybreaks{
\begin{enumerate}
\item
$\displaystyle{T^2_{L}\cong \mathbb{R}^2/L}$, where 
\[L=L_2 =\{n{\bf u}_1 +m{\bf u}_2 =(n+ m a_{1,2}, m a_{2,2})
\,|\,n,m \in\mathbb{Z}\}\]
\[=[\{{\bf u}_1,{\bf u}_2\}]
~\text{is a lattice in} ~\mathbb{R}^2,
\]
\item
$\displaystyle{T^3\cong \mathbb{R}^3/L_3,\,L_3=
\{ n{\bf u}_1 + m{\bf u}_2 + l{\bf u}_3\,|\, n,m,l\in\mathbb{Z}}\}
=[\{{\bf u}_1,{\bf u}_2,{\bf u}_3\}]$,

\item
$\displaystyle{T^4\cong  \mathbb{R}^4/L_4,\,L_4=
\{ n{\bf u}_1 + m{\bf u}_2 + l{\bf u}_3 +k{\bf u}_4
\,|\, n,m,l,k\in \mathbb{Z}\}}=[\{{\bf u}_1,{\bf u}_2,{\bf u}_3,{\bf u}_4\}]$,
\end{enumerate}} 
of two, three and four dimensions. All flat
tori of such dimensions reduce to these cases.
\bigskip

{\bf I. Kroneker's second limit formula and two dimensional tori} 

Since the dual lattice $L^*$ of $L$
is given by
\[
L^*=\left\{ 
\left(n,\frac{m-nA}{B}\right)
\in\mathbb{R}^2\,\mid\, n,m\in \mathbb{Z}
\right\},
\]
non-zero eigenvalues of the Laplacian on 
$\displaystyle{T^2_L}$
are 
\[ 
\left\{
4\pi^2
\left(n^2 + \frac{\left(m-nA\right)^2}{B^2}\right)\,\left|\, 
n,m\in \mathbb{Z}, 
(n,m)\not= (0,0)\right.\right\}.   
\]
The spectral zeta-function $\zeta_{T^2_L}(s)$, 
is 
\allowdisplaybreaks{
\begin{align}
&\frac{1}{\Gamma(s)}
\int_0^{\infty}\left(Z_{T^2_L}\left(t\right)-1\right)t^{s-1}dt
=\zeta_{T^2_L}(s)\notag\\
&= \frac{2}{(4\pi^2)^s}
\sum\limits_{n=1}^{\infty}\sum\limits_{m\in\mathbb{Z}}
\frac{1}{\left(n^2+\left(\frac{1}{B^2}\left(m-nA\right)^2\right)\right)^s}
+\frac{2 B^{2s}}{(4\pi^2)^s}\cdot\zeta(2s)\label{T^2(A,B)}.
\end{align}
}
By using the 
Poisson's summation formula we rewrite the first term as follows:
\allowdisplaybreaks{
\begin{align*}
&\frac{2}{(4\pi^2)^s}
\sum\limits_{n=1}^{\infty}\sum\limits_{m\in\mathbb{Z}}
\frac{1}{\left(n^2+\left(\frac{1}{B^2}
\left(m-nA\right)^2\right)\right)^s}\\
&=\frac{2}{(4\pi^2)^s\cdot\Gamma(s)}
\sum\limits_{n=1}^{\infty}\sum\limits_{m\in\mathbb{Z}}
\frac{1}{n^{2s}}
\int_0^{\infty}
e^{-\left(1+\frac{1}{B^2n^2}\left(m-nA\right)^2\right)x}x^{s-1}dx\\
&=\frac{2}{(4\pi^2)^s\cdot\Gamma(s)}
\sum\limits_{n=1}^{\infty}
\frac{1}{n^{2s}}
\int_0^{\infty}
\sum_{m\in\mathbb{Z}}
\sqrt{\frac{\pi B^2n^2}{x}}
e^{-\frac{(\pi B nm)^2}{x}}e^{-2\pi\sqrt{-1} A n m}
e^{-x}x^{s-1}dx\\
&=\frac{2}{(4\pi^2)^s\cdot\Gamma(s)}
\sum\limits_{n=1}^{\infty}
\frac{1}{n^{2s}}
\int_0^{\infty}
\sum\limits_{m\in\mathbb{Z}, m\not=0}
\sqrt{\frac{\pi B^2n^2}{x}}
e^{-\frac{(\pi B nm)^2}{x}}
e^{-2\pi\sqrt{-1}Anm}
e^{-x}x^{s-1}dx\\
&+\frac{2\sqrt{\pi}B\cdot\Gamma(s-1/2)}{(4\pi^2)^s\cdot\Gamma(s)}
\cdot\zeta(2s-1),
\end{align*}
}
so $\displaystyle{\zeta_{T^2_L}(s)}$ 
is of the following form:
\begin{prop}{\bf (\cite{Mo}, \cite{Be2})}
\label{T^2-zeta decomposition}
\allowdisplaybreaks{
\begin{equation*}
\zeta_{T^2_L}(s)= \mathcal{H}_0(s)
+\frac{2\sqrt{\pi}B\cdot\Gamma(s-1/2)}{(4\pi^2)^s\cdot\Gamma(s)}\cdot
\zeta(2s-1)
+\frac{2 B^{2s}}{(4\pi^2)^s}\cdot\zeta(2s),
\end{equation*}
}
where we put 
\allowdisplaybreaks{
\begin{align*}
\mathcal{H}_0(s)&= \frac{2}{(4\pi^2)^s\cdot\Gamma(s)}
\sum\limits_{n=1}^{\infty} 
\frac{1}{n^{2s}} \int_0^{\infty}
\sum\limits_{m\in\mathbb{Z},\, m\not=0} 
\sqrt{\frac{\pi B^2n^2}{x}}
e^{-\frac{\left(\pi B n m\right)^2}{x}}
e^{-2\pi\sqrt{-1}A n m}
e^{-x}x^{s-1}dx\\
&=\frac{2\sqrt{\pi}B}{(4\pi^2)^s\cdot\Gamma(s)}
\sum\limits_{n=1}^{\infty}  
\frac{e^{-2\pi\sqrt{-1}A n m}}{n^{2s}} 
\int_0^{\infty} 
\sum\limits_{m\in\mathbb{Z},\,m\not=0}
e^{-\frac{\left(\pi B n m\right)^2}{x}}
e^{-x}x^{s-3/2}dx\\
\end{align*}
}
\end{prop}
We know the integrand of $\mathcal{H}_0(s)$ 
satisfies the asymptotics:
\begin{lem}
For any $N\in\mathbb{N}$ 
\begin{equation}\label{infinite-order1}
\sum\limits_{m\in\mathbb{Z},m\not=0}
\sqrt{\frac{\pi B^2n^2}{x}}e^{-\frac{\left(\pi B n m\right)^2}{x}}
e^{-2\pi\sqrt{-1}A n m}
e^{-x}x^{s-1}
=O\left(\frac{x^{N-\mathfrak{Re}(s)-3/2}}{n^{2N-1}}\right).
\end{equation}
\end{lem}
Hence the first term $\mathcal{H}_0(s)$ 
is a holomorphic function of $s$ 
on the whole complex plane, 
and 
\allowdisplaybreaks{
\begin{align*}
&\mathcal{H}_0(s)\\
&=2\sqrt{\pi}B\sum\limits_{n=1}^{\infty}n
\int_0^{\infty}\sum\limits_{m\in\mathbb{Z}\,,\,m\not=0}
e^{-\frac{\left(\pi B n m \right)^2}{x}}e^{-2\pi\sqrt{-1}A n m}
e^{-x}x^{-1/2-1}dx \cdot s +O(s^2) 
\end{align*}
}
Then again by making use of (\ref{Bessel_1/2})
\allowdisplaybreaks{
\begin{align}
\mathcal{H}_0'(0)&=
2\sqrt{\pi}B\sum\limits_{n=1}^{\infty}n e^{-2\pi\sqrt{-1}A n m}
\int_0^{\infty}\sum\limits_{m\in\mathbb{Z}\,,\,m\not=0}
e^{-\frac{\left(\pi B n m\right)^2}{x}}
e^{-x}x^{-1/2-1}dx\notag\\
&=\frac{2}{\sqrt{\pi}}\sum\limits_{n=1}^{\infty}
\sum\limits_{m\in\mathbb{Z},\,m\not=0}
\frac{e^{-2\pi\sqrt{-1}A n m}}{m}\int_0^{\infty}
e^{-\frac{\left(\pi B n m\right)^2}{x}}
e^{-x}x^{-1/2}dx\notag\\
&=2\sum\limits_{n=1}^{\infty}
\sum\limits_{m=1}^{\infty}\frac{1}{m}\left\{
e^{-2\pi n m\left(B-\sqrt{-1}A\right)}
+e^{-2\pi n m \left(B+\sqrt{-1}A\right)}\right\}\notag\\
&=-2\sum\limits_{n=1}^{\infty}
\left\{
\log \left(1-e^{2\pi n\left(B-\sqrt{-1}A\right)}\right)
+\log \left(1-e^{2\pi n\left(B+\sqrt{-1}A\right)}\right)
\right\}.\label{bessel1/2}
\end{align}
}
{}From the facts 
\begin{align*}
&\zeta(-1)=-\frac{1}{12}\\
&\zeta(s)=-\frac{1}{2}-\frac{s}{2}\log 2\pi+ O(s^2)
\end{align*}
and (\ref{Bessel_1/2}) we have
\allowdisplaybreaks{
\begin{align*}
&\zeta_{T^2_L}(s)= -1\\
&+\left\{\frac{\pi B}{3}-2\log B -2\sum\limits_{n=1}^{\infty}
\log\left(1- e^{-2\pi n\left(B-\sqrt{-1}A\right)}\right)+
\log\left(1- e^{-2\pi n\left(B+\sqrt{-1}A\right)}\right)
\right\}\cdot s\\
&+O(s^2).
\end{align*}
}

Then the zeta-regularized determinant 
$Det\,\, \Delta_{T^2_L}$ is given by the formula:
\begin{thm}
\begin{equation}\label{Det T^2}
Det\,\,\Delta_{T^2_L}=
B^2 e^{-\frac{\pi B}{3}}
\prod\limits_{n=1}^{\infty}
\left|
\left(1-e^{-2\pi n(B-\sqrt{-1}A)}\right)
\right|^{4}.
\end{equation}
\end{thm}

\begin{cor} 
{}From the expression {\bf (}\ref{Det T^2}{\bf )} 
we can see easily that
$Det\,\,\Delta_{T^2_L}$ is periodic with respect 
to the parameter
$A$, and when $A=0$ we have both of
\[
\lim\limits_{B\to 0}Det\,\,\Delta_{T^2_L}=0, 
\]
\[\lim\limits_{B\to \infty}Det\,\,\Delta_{T^2_L}=0.
\]
\end{cor}

\smallskip

Now we explain the relation of (\ref{Det T^2})
with the {\it Kroneker's second limit formula}(see \cite{Fo} 
where another explanation is given.).
By using the integral representation of the
modified Bessel function $K_{\alpha}(z)$ (\cite{AAR}):   
\begin{equation*}
K_{\alpha}(z)= \frac{1}{2}\left(\frac{z}{2}\right)^{\alpha}
\int_0^{\infty}e^{-t-\frac{z^2}{4t}}t^{-\alpha-1}dt,\, 
\left|\arg z\right|< \frac{\pi}{4},
\end{equation*}
the function $\mathcal{H}_0(s)$ is expressed as
\begin{equation}
\mathcal{H}_0(s)=
\frac{8B^{s+1/2}}{(4\pi)^s\cdot\Gamma(s)}
\sum\limits_{n=1}^{\infty}
\sum\limits_{m=1}^{\infty}
\cos \left(2\pi\sqrt{-1}A n m\right)
\left(\frac{m}{n}\right)^{s-1/2}
K_{1/2-s}(2\pi B n m).
\end{equation}
Then from this 
we have a functional relation of the function 
$\mathcal{H}_0(s)$:
\begin{prop}\label{H_0-functional-relation}
\begin{equation}
\mathcal{H}_0(1-s)=
\frac{(4\pi)^{2s-1}\cdot\Gamma(s)}{\Gamma(1-s)\cdot B^{2s-1}}\mathcal{H}_0(s),
\end{equation}
especially
\begin{equation}\label{H(1)and H_0'(0)}
\mathcal{H}_0(1)=\frac{B}{4\pi}\mathcal{H}_0'(0).
\end{equation}
\end{prop}

This relation (\ref{H_0-functional-relation}) together
with the functional relation of Riemann $\zeta$-function
gives us a very simple
functional relation of the spectral zeta-function $\zeta_{T^2_L}(s)$ 
for two dimensional flat torus $T^2_L$.
\begin{cor}\label{T2functionalrelation}
\[\Gamma(1-s)\cdot\zeta_{T^2_L}(1-s)=
\left(\frac{4\pi}{B}\right)^{2s-1}\cdot\Gamma(s)\cdot
\zeta_{T^2_L}(s).
\]
\end{cor}

{}From the formula (\ref{T^2-zeta decomposition}) we can easily
see that the function $\zeta_{T^2_L}(s)$ has (only)a pole 
of order one at $s=1$, which comes from that of the second term  
in (\ref{T^2-zeta decomposition}) and the
{\it Kroneker's second limit formula} gives the constant term
at this pole, that is, by the above 
relation (\ref{H(1)and H_0'(0)})
of the first term $\mathcal{H}_0(s)$ we have 
\begin{prop}{\bf (}{\it Kroneker's second limit formula}{\bf )}
\allowdisplaybreaks{
\begin{align*}
&\lim\limits_{s\to 1}\left\{
\zeta_{T^2_L}(s)-\frac{1}{2s-2}\right\}\\
&=\mathcal{H}_0(1)+
\lim\limits_{s\to 1}
\frac{2\sqrt{\pi}B\cdot\Gamma(s-1/2)}{(4\pi)^{s}\cdot\Gamma(s)}
\left\{\zeta(2s-1)-\frac{1}{2s-2}\right\}
+\frac{2B^2}{4\pi^2}\cdot\zeta(2)\\
&=\frac{B}{4\pi}\mathcal{H}_0'(0)
+\frac{B}{2}{\bf C_e}
+\frac{2B^2}{4\pi^2}\cdot\zeta(2).
\end{align*}}
\end{prop}
This gives
\begin{cor}
\[\log Det\,\, \Delta_{T^2}- \frac{4\pi}{B}\lim\limits_{s\to 1}\left\{
\zeta_{T^2_L}(s)-\frac{1}{2s-2}\right\}
=2\log B -2\pi{\bf C_e}
-B\left\{\frac{\pi}{3}+\frac{2}{\pi}\right\}.
\]
\end{cor}

{\bf II. Three dimensional flat torus}.
\bigskip

Let ${\mathfrak A} =\begin{pmatrix}
1 & a_{1,2} & a_{1,3}\\
0 & a_{2,2} & a_{2,3}\\
0 &   0     & a_{3,3} 
\end{pmatrix}$ and ${\mathfrak G}= {^t}{\mathfrak A}^{-1}=(g_{i,j})$,  
then the dual lattice of $L_3$ is generated by
the basis $\{ {\bf u}_1^*,\,{\bf u}_2^*,\, {\bf u}_3^*\}$, where
${\bf u}_j^*=\sum_i g_{i,j}{\bf e}_i^*$ and we know
\allowdisplaybreaks{
\begin{align*}
&Spec(\Delta_{T^3_L})\\
&=\left\{4\pi^2\left(\left(lg_{3,3}+mg_{3,2}+ng_{3,1}\right)^2
+\left(mg_{2,2}+ng_{2,1}\right)^2+
\left(ng_{1,1}\right)^2\right)\,\left|\right.\, n,m,l \in\mathbb{Z}\right\}
\end{align*}
}
Then the spectral zeta-function $\zeta_{T^3_L}(s)$ is 
\allowdisplaybreaks{
\begin{align*}
&\frac{1}{(4\pi^2)^s}
\sum_{{\scriptstyle n,m,l\in\mathbb{Z}}\atop{n^2+m^2+l^2\not=0}} 
\frac{1}{\left(\left(lg_{3,3}+mg_{3,2}+ng_{3,1}\right)^2
+\left(mg_{2,2}+ng_{2,1}\right)^2+
\left(ng_{1,1}\right)^2\right)^s}.
\end{align*}
}
Put $\left(mg_{2,2}+ng_{2,1}\right)^2+
\left(ng_{1,1}\right)^2=I\left(n,m\right)$ and 
as before we express $\zeta_{T^3_L}(s)$ as
\allowdisplaybreaks{
\begin{align*}
&\zeta_{T^3_L}(s)\\
&=\frac{1}{\left(4\pi^2\right)^s}\frac{1}{\Gamma\left(s\right)}
\sum_{{\scriptstyle n,m,\in\mathbb{Z}}\atop{n^2+m^2\not=0}}
\frac{1}{I\left(n,m\right)^s}\sum\limits_{l\in\mathbb{Z}}
\int_0^{\infty}
e^{
\left(-1+ \frac{\left(lg_{3,3}+mg_{3,2}+ng_{3,1}\right)^2}
{I\left(n,m\right)}\right)x
}
x^{s-1}dx\\
&+\frac{2}{\left(2\pi g_{3,3}\right)^{2s}}\cdot\zeta(2s).
\end{align*}
}
Then this equals to the expression: 
\allowdisplaybreaks{
\begin{align*}
&\frac{1}{(4\pi^2)^s}\frac{1}{\Gamma(s)}
\sum_{{\scriptstyle n,m,\in\mathbb{Z}}\atop{n^2+m^2\not=0}}
\frac{1}{I\left(n,m\right)^s}
\int_0^{\infty}
\sum\limits_{l\in\mathbb{Z}}
e^{-\left(l+m\frac{g_{3,2}}{g_{3,3}}+n\frac{g_{3,1}}{g_{3,3}}\right)^2
\frac{g_{3,3}^{\quad 2}}{I\left(n,m\right)}x}e^{-x}x^{s-1}dx\\
&+\frac{2}{\left(2\pi g_{3,3}\right)^{2s}}\cdot\zeta(2s)\\
&=\frac{1}{\left(4\pi^2\right)^s}\frac{1}{\Gamma\left(s\right)}
\sum_{{\scriptstyle n,m,\in\mathbb{Z}}\atop{n^2+m^2\not=0}}
\frac{1}{I\left(n,m\right)^s}
\int_0^{\infty}
\sqrt{\frac{\pi I\left(n,m\right)}{g_{3,3}^{\quad 2}x}}\times\\
&\qquad
\times\sum\limits_{l\in\mathbb{Z}}
e^{-\frac{\pi^2 I\left(n,m\right) l^2}{g_{3,3}^{\,\,\,2}x}} 
e^{2\pi\sqrt{-1}l 
\left(m\frac{g_{3,2}}{g_{3,3}}+n\frac{g_{3,1}}{g_{3,3}}\right)}
e^{-x}x^{s-1}dx
+ \frac{2}{(2\pi g_{3,3})^{2s}}\cdot\zeta(2s)\\
&=\frac{\sqrt{\pi}}{g_{3,3}(4\pi^2)^s}\frac{1}{\Gamma(s)}\times\\
&\times\sum_{{\scriptstyle n,m,\in\mathbb{Z}}\atop{n^2+m^2\not=0}}
\frac{1}{I\left(n,m\right)^{s-1/2}}
\sum_{{\scriptstyle l\in\mathbb{Z}}\atop{l\not=0}}
e^{
2\pi\sqrt{-1}l
\left(
m\frac{g_{3,2}}{g_{3,3}}
+n\frac{g_{3,1}}{g_{3,3}}
\right)}
\int_0^{\infty}
e^{-\frac{\pi^2 I\left(n,m\right) l^2}{g_{3,3}^{\,\,\,2}x}} 
e^{-x}x^{s-3/2}dx\\
&+\frac{\sqrt{\pi}}{g_{3,3}(4\pi^2)^s}
\frac{\Gamma\left(s-1/2\right)}{\Gamma\left(s\right)}
\sum_{{\scriptstyle n,m\in\mathbb{Z}}\atop{n^2+m^2\not=0}}
\frac{1}{I\left(n,m\right)^{s-1/2}}
+\frac{2}{\left(2\pi g_{3,3}\right)^{2s}}\cdot\zeta(2s).
\end{align*}
}
We put this as 
$\mathcal{A}_0(s)+ \mathcal{A}_1(s)+ \mathcal{A}_2(s)$, and calculate
each $\mathcal{A}_i'(0)$.

\begin{enumerate}
\item
\allowdisplaybreaks{
\begin{align*}
&\mathcal{A}_0(0)=0\\
&\mathcal{A}_0'(s)=
-\sum_{{\scriptstyle n,m\in\mathbb{Z}}\atop{n^2+m^2\not=0}}
\log \left(
1-e^{-2\pi\frac{\sqrt{I\left(n,m\right)}}{g_{3,3}}
+ 2\pi\sqrt{-1}\left(m\frac{g_{3,2}}{g_{3,3}}
+n\frac{g_{3,1}}{g_{3,3}}\right)}
\right)\\
&-\sum_{{\scriptstyle n,m\in\mathbb{Z}}\atop{n^2+m^2\not=0}}
\log \left(
1-e^{-2\pi\frac{\sqrt{I\left(n,m\right)}}{g_{3,3}}
-2\pi\sqrt{-1}\left(m\frac{g_{3,2}}{g_{3,3}}
+n\frac{g_{3,1}}{g_{3,3}}\right)}
\right)\\
&=-2\sum\limits_{n=1}^{\infty}
\sum\limits_{m=1}^{\infty}
\log \left(
1-e^{-2\pi\frac{\sqrt{I\left(n,m\right)}}{g_{3,3}}
\pm 2\pi\sqrt{-1}\left(m\frac{g_{3,2}}{g_{3,3}}
+n\frac{g_{3,1}}{g_{3,3}}\right)}
\right)\\
&-2\sum\limits_{n=1}^{\infty}
\sum\limits_{m=1}^{\infty}
\log \left(
1-e^{-2\pi\frac{\sqrt{I\left(n,-m\right)}}{g_{3,3}}
\pm 2\pi\sqrt{-1}\left(-m\frac{g_{3,2}}{g_{3,3}}
+n\frac{g_{3,1}}{g_{3,3}}\right)}\right)\\
&-2\sum\limits_{n=1}^{\infty}
\log 
\left(1-e^{-2\pi n\left\{
\frac{\sqrt{g_{2,1}^2+g_{1,1}^2}}{g_{3,3}}
\pm \sqrt{-1}\frac{g_{3,1}}{g_{3,3}}
\right\}}
\right)\\
&-2\sum\limits_{m=1}^{\infty}
\log \left(
1-e^{-2\pi m\left\{\frac{g_{2,2}}{g_{3,3}}
\pm \sqrt{-1}\frac{g_{3,2}}{g_{3,3}}
\right\}}
\right)
\end{align*}
}
\item
\allowdisplaybreaks{
\begin{align*}
&\mathcal{A}_1(s)
=\frac{\sqrt{\pi}}{g_{3,3}(4\pi^2)^s}\frac{\Gamma(s-1/2)}{\Gamma(s)}
\sum_{{\scriptstyle n,m\in\mathbb{Z}}\atop{n^2+m^2\not=0}}
\frac{1}{I(n,m)^{s-1/2}}\\
&=\frac{1}{2\sqrt{\pi}g_{1,1}^{\quad 2s-1}g_{3,3}}
\frac{\Gamma(s+1/2)}{\Gamma(s+1)}
\frac{s}{s-1/2}
\cdot\zeta_{T^2_L}(s-1/2),
\end{align*}
}
where  $\zeta_{T^2_L}(s)$ is the spectral 
zeta-function of two dimensional
torus $T^2_L$ corresponding to
the lattice 
$L=L_{A,B}$, $\displaystyle{A=-\frac{g_{2,1}}{g_{2,2}}}$, 
$\displaystyle{B=\frac{g_{1,1}}{g_{2,2}}}$.  Then
\allowdisplaybreaks{
\begin{align*}
&\mathcal{A}_1(0)
=0\qquad\qquad\qquad\qquad\qquad\qquad\qquad\qquad\\
&\mathcal{A}_1'(0)
= -\frac{g_{1,1}}{g_{3,3}}\cdot\zeta_{T^2_L}(-1/2).
\end{align*}
}
Next we express the value $\zeta_{T^2_L}(-1/2)$
in terms of a modified Bessel function $K_{\alpha}(z)$ for $\alpha = 1$:
\allowdisplaybreaks{
\[
K_{\alpha}(z)= \frac{1}{2}\left(\frac{z}{2}\right)^{\alpha}
\int_0^{\infty}e^{-t-\frac{z^2}{4t}}t^{-\alpha-1}dt,\, 
\left|\arg z\right|< \frac{\pi}{4}.
\] 
}
So we return to the expression
(\ref{T^2(A,B)}):
\allowdisplaybreaks{
\begin{equation*}
\zeta_{T^2_L}(s)= \mathcal{H}_0(s)
+\frac{2\sqrt{\pi}B\cdot\Gamma(s-1/2)}{(4\pi^2)^s\cdot\Gamma(s)}\cdot\zeta(2s-1)
+\frac{2 B^{2s}}{(4\pi^2)^s}\cdot\zeta(2s).
\end{equation*}
}
Then
\allowdisplaybreaks{
\begin{align*}
&\mathcal{H}_0(-1/2)\\
&=-2\pi B
\sum\limits_{n=1}^{\infty}
n^2\int_0^{\infty}\sum\limits_{m\in\mathbb{Z},\,m\not=0}
e^{-\frac{(\pi B n m)^2}{x}}
e^{-2\pi\sqrt{-1}A n m}
e^{-x}
x^{-2}dx\\
&=-8\sum\limits_{n=1}^{\infty}
\sum\limits_{m=1}^{\infty}(\frac{n}{m})
K_{1}(2\pi B n m)\cos(2\pi A n m),\\
&\frac{2\sqrt{\pi}B\cdot\Gamma(s-1/2)}{(4\pi^2)^s\cdot\Gamma(s)}
\cdot\zeta(2s-1)_{|s=-1/2}=
-\frac{2B}{\pi}\cdot\zeta(3),\\
&\frac{2 B^{2s}}{(4\pi^2)^s}\cdot\zeta(2s)_{|s=-1/2}=
-\frac{\pi}{3B}.
\end{align*}
}
Hence
\allowdisplaybreaks{
\[\mathcal{A}_1'(0)=\frac{g_{1,1}}{g_{3,3}}
\left\{8\sum\limits_{n=1}^{\infty}
\sum\limits_{m=1}^{\infty}\frac{n}{m}
K_{1}(2\pi B n m)\cos\left(2\pi A n m\right)
+\frac{2B}{\pi}\cdot\zeta(3)+\frac{\pi}{3B}\right\}.
\]
}
\item
\allowdisplaybreaks{\begin{align*}
&\mathcal{A}_2(0)=-1
\qquad\qquad\qquad\qquad\qquad
\qquad\qquad\qquad\qquad\qquad\qquad\qquad\\
&\mathcal{A}_2'(0)=2\log g_{3,3}.
\end{align*}}
\end{enumerate}

Finally, for the lattice $L_3=[\{{\bf u}_i\}_{i=1}^{3}]$
we have an expression
of the zeta-regularized determinant
\begin{thm}
\allowdisplaybreaks{\begin{align*}
&Det \,\, \Delta_{T^3_{L}}\\
&=\prod\limits_{n=1}^{\infty}
\prod\limits_{m=1}^{\infty}
\left|
1-e^{-2\pi
\left\{
\frac{\sqrt{I(n,m)}}{g_{3,3}}
+\sqrt{-1}\frac{m g_{3,2}+ n g_{3,1}}{g_{3,3}}
\right\}}\right|^4\times\\
&\times\prod\limits_{n=1}^{\infty}
\prod\limits_{m=1}^{\infty}
\left|
1-e^{-2\pi
\left\{
\frac{\sqrt{I(n,-m)}}{g_{3,3}}
+\sqrt{-1}\frac{-m g_{3,2}+ng_{3,1}}{g_{3,3}}
\right\}}
\right|^4\times\\
&\times\prod\limits_{n=1}^{\infty}
\left|
1-e^{-2\pi n
\left\{\frac{\sqrt{g_{2,1}^2+g_{1,1}^2}}{g_{3,3}}
+\sqrt{-1}\frac{g_{3,1}}{g_{3,3}}
\right\}}
\right|^4\cdot\prod\limits_{m=1}^{\infty}
\left|
1-e^{-2\pi m
\left\{\frac{g_{2,2}}{g_{3,3}}
+\sqrt{-1}\frac{g_{3,2}}{g_{3,3}}
\right\}}
\right|^4\times\\
&\times
e^{-\frac{g_{1,1}}{g_{3,3}}\left\{
\sum\limits_{n=1}^{\infty}
\sum\limits_{m=1}^{\infty}\left(\frac{n}{m}\right)
K_{1}\left(2\pi\frac{g_{1,1}}{g_{2,2}} n m\right)
\cos\left(2\pi\frac{g_{2,1}}{g_{2,2}} n m\right)\right\}}
\cdot
e^{-\frac{2}{\pi}\frac{g_{1,1}^{\quad 2}}{g_{3,3}g_{2,2}}\cdot\zeta(3)}
\cdot
e^{-\frac{\pi}{3}\frac{g_{2,2}}{{g_{3,3}}}}
\times\\
&\times\left(\frac{1}{g_{3,3}}\right)^2.
\end{align*}}
\end{thm}


\begin{rem}
As a special case of $T^3$, let assume $a_{1,3}=a_{2,3}=0$. Then 
the formula of $Det\,\,\Delta_{T^3}$ for this case
coincides with the formula (\ref{general_product_formula})
for $T^2\times S^1$.
\end{rem}

\bigskip

{\bf III. Four dimensional flat torus}

Here we only state a formula for a four dimensional torus
and a special case of them.

\bigskip

Let ${\mathfrak A}=\begin{pmatrix}
1 & a_{1,2} & a_{1,3}& a_{1,4}\\
0 & a_{2,2} & a_{2,3}& a_{2,4}\\
0 &   0     & a_{3,3}& a_{3,4}\\
0 &   0     &    0   & a_{4,4} 
\end{pmatrix}$ and ${\bf u}_j=\sum\limits_i a_{i,j}{\bf e}_i$ 
as explained in the beginning of this section
and ${\mathfrak G}= {^t}{\mathfrak A}^{-1}=(g_{i,j})$.  

Let $L=L({\mathfrak A})$ be the lattice generated 
by $\{{\bf u}_1,{\bf u}_2,{\bf u}_3,{\bf u}_4\}$,
then the dual lattice of $L({\mathfrak A})$ is generated by
the basis $\{ {\bf u}_1^*,\,{\bf u}_2^*,\, {\bf u}_3^*,\, {\bf u}_4^*\}$, where
${\bf u}_j^*=\sum_i g_{i,j}{\bf e}_i^*$.

Put 
\[\left(lg_{3,3}+mg_{3,2}+ng_{3,1}\right)^2
+\left(mg_{2,2}+ng_{2,1}\right)^2+
\left(ng_{1,1}\right)^2=I(n,m,l),
\]
and
\[
\frac{ng_{4,1}+mg_{4,2}+lg_{4,3}}{g_{4,4}}=\alpha(n,m,l)
\]
then the spectral zeta-function $\zeta_{T^4_{L({\mathfrak A})}}(s)$ 
for this case is written as
\allowdisplaybreaks{
\begin{align*}
&\zeta_{T^4_{L({\mathfrak A})}}(s)\\
&=\frac{1}{\left(4\pi^2\right)^s}
\sum_{{\scriptstyle n,m,l\in\mathbb{Z}}\atop{n^2+m^2+l^2\not=0}} 
\frac{1}
{\left(I\left(n,m,l\right)
+\left(ng_{4,1}+mg_{4,2}+lg_{4,3}+kg_{4,4}\right)^2\right)^s}\\
&=\frac{1}{(4\pi^2)^s}\frac{1}{\Gamma(s)}
\sum_{{\scriptstyle n,m,l\in\mathbb{Z}}\atop{n^2+m^2+l^2\not=0}}
\frac{1}{I\left(n,m,l\right)^s}
\int_0^{\infty}\sum\limits_{k\in\mathbb{Z}}
e^{-\left(k+\alpha\left(n,m,l\right)\right)^2
\frac{g_{4,4}^{\quad 2}}{I\left(n,m,l\right)}x}
e^{-x}
x^{s-1}dx\\
&+\frac{2}{\left(2\pi g_{4,4}\right)^{2s}}\cdot\zeta(2s)\\
&=\frac{\sqrt{\pi}}{g_{4,4}\left(4\pi^2\right)^s}
\frac{1}{\Gamma(s)}\times\\
&\times\sum_{{\scriptstyle n,m,l\in\mathbb{Z}}\atop{n^2+m^2+l^2\not=0}}
\frac{1}{I(n,m,l)^{s-1/2}}
\int_0^{\infty}\sum\limits_{k\in\mathbb{Z},k\not=0}
e^{-\frac{\pi^2 k^2}{xg_{4,4}^{\quad 2}}I\left(n,m,l\right)}
e^{2\pi\sqrt{-1}\alpha\left(n,m,l\right) k}
e^{-x}x^{s-3/2}
dx\\
&+\frac{\sqrt{\pi}}{g_{4,4}\left(2\pi\right)^{2s}}
\frac{\Gamma(s-1/2)}{\Gamma(s)}
\sum_{{\scriptstyle n,m,l\in\mathbb{Z}}\atop{n^2+m^2+l^2\not=0}}
\frac{1}{I\left(n,m,l\right)^{s-1/2}}
+\frac{2}{\left(2\pi g_{4,4}\right)^{2s}}\cdot\zeta(2s).
\end{align*}
}
We put this as 
$\mathcal{B}_0(s)+\mathcal{B}_1(s)+\mathcal{B}_2(s)$ 
corresponding to each term.

Note that
\[
\mathcal{B}_1(s)
=\frac{1}{2\sqrt{\pi}g_{4,4}}\frac{\Gamma(s+1/2)}{\Gamma(s+1)}
\cdot\frac{s}{s-1/2}\cdot \zeta_{T^3}(s-1/2),
\]
and at $s=-1/2$, $\zeta_{T^3}(s)$ is holomorphic. 

\begin{thm}
$$
Det\,\,\Delta_{T^4_{L({\mathfrak A})}}=\prod\limits_{i=0}^2
e^{-\mathcal{B}_i'(0)},
$$ 
where each $\mathcal{B}_i'(0)$ is given as follows:
\allowdisplaybreaks{
\begin{align*}
&\mathcal{B}_0(0)=0,\\
&\mathcal{B}_0'(0)\\
&\qquad=-\sum_{{\scriptstyle n,m,l\in\mathbb{Z}}\atop{n^2+m^2+l^2\not=0}}
\log \left(1-e^{-2\pi \left\{\frac{\sqrt{I\left(n,m,l\right)}}{g_{4,4}}
+\sqrt{-1}\alpha\left(n,m,l\right)\right\}}\right)\\
&\qquad\qquad \qquad \qquad 
-\sum_{{\scriptstyle n,m,l\in\mathbb{Z}}\atop{n^2+m^2+l^2\not=0}}
\log \left(1-e^{-2\pi \left\{\frac{\sqrt{I\left(n,m,l\right)}}{g_{4,4}}
-\sqrt{-1}\alpha\left(n,m,l\right)\right\}}\right),\\
&\mathcal{B}_1(0)=0,\\
&\mathcal{B}_1'(0)=-\frac{1}{g_{4,4}}\cdot\zeta_{T^3}(-1/2)
=-\frac{1}{g_{4,4}}
\left\{
{\mathcal A}_0\left(-1/2\right)
+{\mathcal A}_1\left(-1/2\right)
+{\mathcal A}_2\left(-1/2\right)
\right\}\\
&=\frac{1}{g_{4,4}}
\left\{\sum_{{\scriptstyle n,m,\in\mathbb{Z}}\atop{n^2+m^2\not=0}}
\sum\limits_{l=1}^{\infty}
4\frac{\sqrt{I\left(n,m\right)}}{l}
\cos \left(2\pi l \frac{mg_{3,2}+ng_{3,1}}{g_{3,3}}\right)
\cdot
K_1\left(\frac{2\pi}{g_{3,3}l}\sqrt{I\left(n,m\right)}\right)\right.\\
&\left.
\qquad\qquad
+8\sqrt{\pi}\frac{g_{1,1}^2}{g_{3,3}}
\sum\limits_{n=1}^{\infty}
\sum\limits_{m=1}^{\infty}
\frac{n^2}{m}\cos \left(2\pi\sqrt{-1}\frac{g_{2,1}}{g_{2,2}}n m\right)
K_1\left(2\pi \frac{g_{1,1}}{g_{2,2}}n m \right)\right.\\
&\left.\qquad\qquad\qquad\qquad\qquad\qquad\qquad\qquad
+\frac{3g_{1,1}^2}{4\pi^3 g_{3,3}}\cdot\zeta(4)
+\frac{2g_{2,2}^{\quad 2}}{\pi g_{3,3}}\cdot\zeta(3)
+\frac{\pi}{3}{g_{3,3}}\right\},\\
&\mathcal{B}_2(0)=-1,\\
&\mathcal{B}_2'(0)=2\log g_{4,4}.
\end{align*}
}
\end{thm}

The formula in
Theorem \ref{Det M*N e-e} and similar one give us several formulas
of the zeta-regularized determinant for flat tori defined by 
matrices of the form 
${\mathfrak A}=\begin{pmatrix}
1 & a_{1,2} & a_{1,3}& 0\\
0 & a_{2,2} & a_{2,3}& 0\\
0 &   0     & a_{3,3}& 0\\
0 &   0     &    0   &a_{4,4} 
\end{pmatrix}$ or ${\mathfrak A}=\begin{pmatrix}
1 & a_{1,2} &       0&       0\\
0 & a_{2,2} &       0&       0\\
0 &   0     & a_{3,3}& a_{3,4}\\
0 &   0     &    0   & a_{4,4} 
\end{pmatrix}$. 

As an application of the formula (\ref{Det M*N e-e})
we state a formula for a torus defined by the latter one, that is,
the torus is a direct product of two 2-dimensional tori 
as Riemannian manifold.  

Let 
${\mathfrak A}=\begin{pmatrix}
1 & a_{1,2} &       0&       0\\
0 & a_{2,2} &       0&       0\\
0 &   0     & a_{3,3}& a_{3,4}\\
0 &   0     &    0   & a_{4,4} 
\end{pmatrix}$ and 
$L_{\mathfrak A}$ the lattice in $\mathbb{R}^4$ generated
by $\{{\bf u}_1={\bf e}_1, {\bf u}_2= a_{1,2}{\bf e}_1+
a_{2,2}{\bf e}_2, {\bf u}_3=a_{3,3}{\bf e}_3,
{\bf u}_4=a_{3,4}{\bf e}_3+a_{4,4}{\bf e}_4\}$. Then 
the torus $\mathbb{R}^4/L({\mathfrak A})$ is a direct product of two tori
${\bf M}\times {\bf N}$, where
${\bf M}=\mathbb{R}^2/L_{\bf M},~L_{\bf M}=[\{{\bf u}_1,{\bf u}_2\}]$ and
${\bf N}=\mathbb{R}^2/L_{\bf N},~L_{\bf N}=[\{{\bf u}_3,{\bf u}_3\}]$.
Since each heat kernel asymptotics for flat tori
vanishes except the first one = ${\bf b}_0$ = {\it volume of the torus}, 
we can rewrite the formula (\ref{Det M*N e-e}) for this case as
\begin{cor}
\begin{equation*}
Det\,\,\Delta_{T^4_{L({\mathfrak A})}}
=Det\,\,\Delta_{\bf M}
\cdot\prod\limits_{\lambda\in L^*_{\bf M}, \lambda\not=0}
e^
{-\int_0^{\infty}
\left( 
K_{\bf N}\left(\frac{t}{\Vert\lambda\Vert^2}\right)
-\frac{\Vert\lambda\Vert^2 {\bf b}_0}{4\pi t}
\right)
e^{-t}t^{-1}dt
}\cdot
e^
{-\frac{{\bf b}_0}{4\pi}\{\lim\limits_{s\to 0}
\Gamma(s-1)\cdot\zeta_{\bf M}(s-1)\}},
\end{equation*}
where 
\begin{align*}
&Det\,\,\Delta_{\bf M}=
B^2 e^{-\frac{\pi B}{3}}
\prod\limits_{n=1}^{\infty}
\left|
\left(1-e^{-2\pi n\left(B-\sqrt{-1}A\right)}\right)
\right|^{4},
\end{align*}
\allowdisplaybreaks{
\begin{align*}
&\prod\limits_{\lambda\in L^*_{\bf M},\lambda\not=0}
e^
{-\int_0^{\infty}
\left( 
K_{\bf N}\left(
\frac{t}{\Vert\lambda\Vert^2}\right)
-\frac{\Vert\lambda\Vert^2 {\bf b}_0}{4\pi t}
\right)
e^{-t}t^{-1}dt}
=\prod\limits_{\lambda\in L^*_{\bf M},\lambda\not=0}
\prod\limits_{\gamma\in L_{\bf N},\gamma\not=0}
e^{-\frac{\Vert\lambda\Vert {\bf b}_0}{\pi\Vert\gamma\Vert}
K_1(\Vert\lambda\Vert \Vert\gamma\Vert)},
\end{align*}}
and
\allowdisplaybreaks{
\begin{align*}
&\lim\limits_{s\to 0}\Gamma(s-1)\cdot\zeta_{\bf M}(s-1)\\
&=\lim\limits_{s\to -1}
\left\{
\frac{2}{(4\pi^2)^s}
\sum\limits_{n=1}^{\infty}
\frac{1}{n^{2s}}
\int_0^{\infty}
\sum\limits_{m\in\mathbb{Z},m\not=0}
\sqrt{\frac{\pi B^2n^2}{x}}e^{-\frac{(\pi B n m)^2}{x}}
e^{-2\pi\sqrt{-1} A n m}
e^{-x}x^{s-1}dx\right.\\
&\left.\qquad\qquad\qquad\qquad
+\frac{2\sqrt{\pi}B}{(4\pi^2)^s}\cdot
\Gamma(s-1/2)\cdot\zeta(2s-1)\right.\\
&\left.\qquad\qquad\qquad\qquad\qquad\qquad\qquad\qquad
+\frac{\Gamma(s)\cdot B^{2s}}{(4\pi^2)^s}\cdot\zeta(2s)\right\}\\
&=\frac{32\pi}{\sqrt{B}}\sum\limits_{n=1}^{\infty}
\sum\limits_{m=1}^{\infty}
\left(\frac{n}{m}\right)^{3/2}
\cos (2\pi Anm)K_{3/2}(2\pi B n m)
+\frac{8\pi}{B}\cdot\zeta(4)+\frac{4}{B^2}\cdot\zeta(3)\\
&=\left(\frac{4\pi}{B}\right)^3\Gamma(2)\zeta_{T^2}(2)=4\pi B
\sum_{{\scriptstyle n,m\in\mathbb{Z}}\atop{n^2+m^2\not=0}}
\frac{1}{\left(\left(Bn\right)^2+\left(m-nA\right)^2\right)^2},
\end{align*}}
here $B=\frac{g_{1,1}}{g_{2,2}}$, $A=-\frac{g_{2,1}}{g_{2,2}}$,
${\bf b}_0$ = {\it volume of} ${\bf N}$ = $a_{3,3}a_{4,4}$.
\end{cor}

Note that the second term 
is obtained by making use of the Jacobi identity:
\begin{align*}
&\sum\limits_{\mu\in L^*_{\bf N}}e^{-t\Vert\mu\Vert^2}
=\frac{{\bf b}_0}{4\pi t}\sum\limits_{\gamma\in L_{\bf N}}
e^{-\frac{\Vert\gamma\Vert^2}{4t}}.
\end{align*} 

Finally, we note that in the most special case, that is, 
let the matrix ${\mathfrak A}$ be the identity
matrix, then 
we have a formula for the spectral zeta-function $\zeta_{T^4}(s)$ 
\begin{equation}
(4\pi^2)^{-s}\zeta_{T^4}(s)=8(1-2^{2-2^s})\zeta(s)\zeta(s-1)
\end{equation}
(and in each dimension we have similar formulas).
So by this formula we have simply
\begin{cor}
Let the torus $T^4$ be defined by the lattice 
$\left\{{\bf e}_1, {\bf e}_2, {\bf e}_3, {\bf e}_4\right\}$, then
the zeta-regularized determinant $Det\,\, \Delta_{T^4}$
is given explicitly in the form
\begin{align*}
&\log Det\,\, \Delta_{T^4}
=2^4 \left(\log 2\pi +2(\log 2)^2\right)\zeta(0)\zeta(-1)
-2^3\left(\zeta'(0)\zeta(-1)+\zeta(0)\zeta'(-1)\right)\\
\end{align*}
\end{cor}
It is possible to simplify this formula by using several formulas
of Riemann  zeta-function $\zeta(s)$ and to compare with our formula.

\begin{rem}
Similar to the case of $\zeta_{T^2_L}(s)$, 
the function $\zeta_{T^3_L}(s)$ (respectively
$\zeta_{T^4_L}(s)$) has only a pole at $s=3/2$ (resp. $s=2$)
of order one coming from the second term $\mathcal{A}_1(s)$
(respectively $\mathcal{B}_1(s)$) and  
the term $\mathcal{A}_0(s)$ (resp. $\mathcal{B}_0(s)$) 
will correspond to the term
$\mathcal{H}_0(s)$ in the two dimensional cases and there are
similar functional relations like (\ref{H_0-functional-relation}) 
also in these cases which are derived from the Jacobi identity. 
\end{rem}


\begin{align*}
\end{align*}

%
%
%
%
%
%
%




\end{document}